\documentclass[12pt,a4paper]{article}
\newtheorem{thm}{Theorem}

\newtheorem{prop}{Proposition}
\newtheorem{cor}{Corollary}
\newcommand{\dimo}{\noindent \emph{Proof. }}
\newcommand{\qed}{\\ \rightline{$\Box$}\\}

\usepackage{latexsym}
\usepackage{amsfonts}
\usepackage{amssymb}
\usepackage{amsmath}
\usepackage{amscd}
\usepackage[dvips]{graphics}
\usepackage[dvips]{epsfig}
\usepackage{graphicx}
\usepackage[english]{babel}

\begin{document}

\title{Cyclic generalizations of two\\ hyperbolic icosahedral manifolds}

\author{P.~Cristofori - T.~Kozlovskaya - A.~Vesnin}




\maketitle

\begin{abstract}
We discuss two families of closed orientable three-dimensional manifolds which arise as cyclic generalizations of two hyperbolic icosahedral manifolds listed by Everitt. Everitt's manifolds are cyclic coverings of the lens space $L_{3,1}$ branched over some 2-component links. We present results on covering properties, fundamental groups, and hyperbolic volumes of the manifolds belonging to these families. 
\end{abstract}

\vskip 30pt
{\bf keywords: } 3-manifold,\  cyclic branched covering,\ lens space,\ links in manifolds

\vskip 20pt


\section*{Introduction}
\label{sec0}

Various examples of three-dimensional spherical, Euclidean, or hyperbolic manifolds arise from pairwise isometrical identifications of faces of convex regular polyhedra in corresponding 3-spaces: $\mathbb S^3$, $\mathbb E^3$, or $\mathbb H^3$. The most famous examples are the spherical and hyperbolic dodecahedral manifolds constructed by Weber and Seifert in 1933 \cite{Weber-Seifert}. The whole set of such examples for every  spherical, Euclidean, or hyperbolic convex regular polyhedron was listed by Everitt \cite{Everitt}.  The list contains eight manifolds $M_{15}, \ldots , M_{22}$, arising from a regular hyperbolic dodecahedron with dihedral angled $2\pi / 5$, and six manifolds $M_{23}, \ldots , M_{28}$, arising from a regular hyperbolic icosahedron with dihedral angles $2 \pi / 3$.  It can be checked directly from the gluing schemata that $M_{15}$ is the Weber~-- Seifert manifold from \cite{Weber-Seifert}, and $M_{23}$ is the Fibonacci manifold from \cite{HKM1}, uniformized by the Fibonacci group $F(2,10)$.  Both manifolds have cyclic symmetries induced by symmetries of the polyhedra, such that $M_{15}$ is an $5$-fold cyclic covering of the 3-sphere $S^{3}$, branched over the Whitehead link, and $M_{23}$ is the $5$-fold cyclic covering of $S^{3}$, branched over the figure-eight knot. Cyclic generalizations of these manifolds were constructed in  \cite{HKM2} and \cite{HKM1}: the $n$-fold strongly-cyclic coverings of $S^{3}$, branched over the Whitehead link, and the $n$-fold cyclic coverings of $S^{3}$, branched over the figure-eight knot, respectively. Explicit  formulae for hyperbolic volumes of manifolds of these two classes  are given in \cite{MedVes2} and \cite{MedVes1}.

It was observed by Cavicchioli, Spaggiari and Telloni \cite{Cavicchioli} that the manifolds $M_{24}$ and $M_{25}$, arising from the $2 \pi / 3$-icosahedron, are $3$-fold cyclic branched coverings of the lens space $L_{3,1}$ branched over some 2-component links. In the present paper we will consider two families of 3-manifolds which are cyclic generalizations of $M_{24}$ and $M_{25}$.

One family of manifolds, namely $M_{24}(n)$, $n \geqslant 1$, is a generalization of the manifold $M_{24}$ from~\cite{Everitt}. The manifolds $M_{24}(n)$, where $M_{24}(3) = M_{24}$, were independently constructed by Cavicchioli, Spaggiari and Telloni~\cite{Cavicchioli2} for $n \geqslant 3$ and by Kozlovskaya~\cite{Kozlovskaya1, Kozlovskaya2} for $n\geqslant 2$. In both cases $M_{24}(n)$ was defined via pairwise identifications of the faces of a 3-complex. These manifolds are natural generalizations of $M_{24}$ in the following sense: $M_{24}(n)$, $n > 1$, is an $n$-fold strongly-cyclic branched covering of the lens space $L_{3,1}$, branched over the same link as $M_{24}$. The proof of this fact, presented in~\cite{Cavicchioli, Cavicchioli2}, is based on results by Stevens \cite{Stevens} and by Osborne and Stevens~\cite{Osborn-Stevens}. In Theorem~\ref{theorem:m24n-covering} we give a purely topological proof of this fact: we consider a Heegaard diagram for the quotient space of $M_{24}(n)$  by its cyclic symmetry and we reduce it to the standard genus one Heegaard diagram  for $L_{3,1}$. The manifold $M_{24}(3)$ is hyperbolic since, by construction, it is obtained from the hyperbolic $2\pi / 3$--icosahedron gluing its faces by isometries.  It is stated in \cite{Cavicchioli2}, without proof, that for $n  > 3$ the manifolds $M_{24}(n)$ are hyperbolic and a formula for their volumes is given. However the formula turns out to be wrong even for $n=4$. In fact, in Proposition~\ref{prop1.4} we will present the correct values of volumes for the initial list of manifolds $M_{24}(n)$, calculated by the computer program \emph{Recognizer}\footnote{Three-manifold Recognizer is a computer program developed by the research group of S.~Matveev in the Department of Computer Topology and Algebra of Chelyabinsk State University, available on the webpage http://www.matlas.math.csu.ru}. In \cite{Ves-Koz} the construction of $M_{24}(n)$ was generalized to obtain $n$-fold cyclic coverings of lens spaces $L(p,q)$, branched over two components.

Another family of manifolds, namely $M_{25}(n)$, $n \geqslant 1$, is a generalization of the manifold $M_{25}$ from \cite{Everitt}. As well as $M_{24}$, $M_{25}$ can be constructed from the $2\pi/3$--icosahedron. It is written in  \cite[p.~931]{Cavicchioli2} that ``one can construct an analogous $\mathbb Z_{n}$- symmetric descriprion for $M_{25}(n)$ in the same manner'' as $M_{24}(n)$, with $M_{25}(3) = M_{25}$. Unfortunately,  \cite{Cavicchioli2} doesn't contain an explicit description of $M_{25}(n)$ for arbitrary $n$, but only the presentation for the fundamental group induced by face pairings. In this paper we reconstruct $M_{25}(n)$ by means of the face pairing corresponding to the presentation for $\pi_{1}(M_{25}(n))$ announced in \cite{Cavicchioli2}. We prove that, for $n$ even, the manifolds $M_{25}(n)$ have covering properties, fundamental group presentations, and hyperbolic volumes different from stated in Theorem 3.1 of  \cite{Cavicchioli2}, see   Theorem~\ref{theorem:m25n-covering}, Proposition~\ref{proposition-m25n-presentation}. 

We recall that cyclic branched coverings of the 3-sphere branched over knots or links are intensively investigated from various points of views (see, for example a survey for the case of 2-bridge knots and links in \cite{Mulazzani-Vesnin}). About existence and uniqueness for strongly-cyclic branched coverings of 3-manifolds branched over knots, see \cite{Cristofori-Mulazzani-Vesnin}. Manifolds $M_{24}(n)$ and $M_{25}(n)$ seem to be very interesting as a starting point to understand polyhedral constructions and Heegaard diagrams for cyclic coverings of 3-manifolds (in particular, lens spaces) branched over links.

\section{The family of manifolds $M_{24}(n)$}
\label{sec1}

\begin{figure}[h]
\begin{center}
\unitlength=0.75mm
\begin{picture}(140,75)(-20,-5)
\thicklines
\qbezier(-20,0)(-20,0)(100,0) \qbezier(0,60)(0,60)(120,60)
\qbezier(20,0)(20,0)(20,40) \qbezier(60,0)(60,0)(60,40)
\qbezier(100,0)(100,0)(100,40)
\qbezier(0,20)(0,20)(0,60) \qbezier(40,20)(40,20)(40,60)
\qbezier(80,20)(80,20)(80,60) \qbezier(120,20)(120,20)(120,60)
\qbezier(20,40)(20,40)(0,20) \qbezier(20,40)(20,40)(0,60)
\qbezier(20,40)(20,40)(40,20)  \qbezier(20,40)(20,40)(40,60)
\qbezier(60,40)(60,40)(40,20) \qbezier(60,40)(60,40)(40,60)
\qbezier(60,40)(60,40)(80,20)  \qbezier(60,40)(60,40)(80,60)
\qbezier(100,40)(100,40)(80,20) \qbezier(100,40)(100,40)(80,60)
\qbezier(100,40)(100,40)(120,20)  \qbezier(100,40)(100,40)(120,60)
\qbezier(0,20)(0,20)(-20,0)
\qbezier(0,20)(0,20)(20,0) \qbezier(40,20)(40,20)(20,0)
\qbezier(40,20)(40,20)(60,0)  \qbezier(80,20)(80,20)(60,0)
\qbezier(80,20)(80,20)(100,0) \qbezier(120,20)(120,20)(100,0)
\put(120,65){\makebox(0,0)[cc]{$P_1$}}
\put(80,65){\makebox(0,0)[cc]{$P_2$}}
\put(40,65){\makebox(0,0)[cc]{$P_3$}}
\put(0,65){\makebox(0,0)[cc]{$P_1$}}
\put(100,47){\makebox(0,0)[cc]{$Q_1$}}
\put(60,47){\makebox(0,0)[cc]{$Q_2$}}
\put(20,47){\makebox(0,0)[cc]{$Q_3$}}
\put(120,13){\makebox(0,0)[cc]{$R_1$}}
\put(80,13){\makebox(0,0)[cc]{$R_2$}}
\put(40,13){\makebox(0,0)[cc]{$R_3$}}
\put(0,13){\makebox(0,0)[cc]{$R_1$}}
\put(100,-5){\makebox(0,0)[cc]{$S_1$}}
\put(60,-5){\makebox(0,0)[cc]{$S_2$}}
\put(20,-5){\makebox(0,0)[cc]{$S_3$}}
\put(-20,-5){\makebox(0,0)[cc]{$S_1$}}
\put(100,55){\makebox(0,0)[cc]{$\bf A_1$}}
\put(60,55){\makebox(0,0)[cc]{$\bf A_2$}}
\put(20,55){\makebox(0,0)[cc]{$\bf A_3$}}
\put(110,40){\makebox(0,0)[cc]{$\bf B_1$}}
\put(90,40){\makebox(0,0)[cc]{$\bf \bar{A}_3$}}
\put(70,40){\makebox(0,0)[cc]{$\bf B_2$}}
\put(50,40){\makebox(0,0)[cc]{$\bf \bar{A}_1$}}
\put(30,40){\makebox(0,0)[cc]{$\bf B_3$}}
\put(10,40){\makebox(0,0)[cc]{$\bf \bar{A}_2$}}
\put(110,20){\makebox(0,0)[cc]{$\bf C_1$}}
\put(90,20){\makebox(0,0)[cc]{$\bf \bar{C}_1$}}
\put(70,20){\makebox(0,0)[cc]{$\bf C_2$}}
\put(50,20){\makebox(0,0)[cc]{$\bf \bar{C}_2$}}
\put(30,20){\makebox(0,0)[cc]{$\bf C_3$}}
\put(10,20){\makebox(0,0)[cc]{$\bf \bar{C}_3$}}
\put(80,5){\makebox(0,0)[cc]{$\bf \bar{B}_1$}}
\put(40,5){\makebox(0,0)[cc]{$\bf \bar{B}_2$}}
\put(0,5){\makebox(0,0)[cc]{$\bf \bar{B}_3$}}
\put(50,70){\makebox(0,0)[cc]{$\bf D$}}
\put(50,-10){\makebox(0,0)[cc]{$\bf \bar{D}$}}
\put(100,65){\makebox(0,0)[cc]{$x_1$}}
\put(102,60){\vector(-1,0){4}} \put(60,65){\makebox(0,0)[cc]{$x_2$}}
\put(62,60){\vector(-1,0){4}} \put(20,65){\makebox(0,0)[cc]{$x_3$}}
\put(22,60){\vector(-1,0){4}} \put(80,-5){\makebox(0,0)[cc]{$x_2$}}
\put(82,0){\vector(-1,0){4}} \put(40,-5){\makebox(0,0)[cc]{$x_3$}}
\put(42,0){\vector(-1,0){4}} \put(0,-5){\makebox(0,0)[cc]{$x_1$}}
\put(2,0){\vector(-1,0){4}} \put(125,40){\makebox(0,0)[cc]{$x_2$}}
\put(120,38){\vector(0,1){4}} \put(83,40){\makebox(0,0)[cc]{$x_3$}}
\put(80,38){\vector(0,1){4}} \put(43,40){\makebox(0,0)[cc]{$x_1$}}
\put(40,38){\vector(0,1){4}} \put(-5,40){\makebox(0,0)[cc]{$x_2$}}
\put(0,38){\vector(0,1){4}} \put(92,52){\makebox(0,0)[cc]{$u$}}
\put(88,52){\vector(1,-1){4}} \put(52,52){\makebox(0,0)[cc]{$u$}}
\put(48,52){\vector(1,-1){4}} \put(12,52){\makebox(0,0)[cc]{$u$}}
\put(8,52){\vector(1,-1){4}} \put(108,52){\makebox(0,0)[cc]{$y_1$}}
\put(112,52){\vector(-1,-1){4}}
\put(68,52){\makebox(0,0)[cc]{$y_2$}} \put(72,52){\vector(-1,-1){4}}
\put(28,52){\makebox(0,0)[cc]{$y_3$}} \put(32,52){\vector(-1,-1){4}}
\put(108,12){\makebox(0,0)[cc]{$y_3$}} \put(108,8){\vector(1,1){4}}
\put(68,12){\makebox(0,0)[cc]{$y_1$}} \put(68,8){\vector(1,1){4}}
\put(28,12){\makebox(0,0)[cc]{$y_2$}} \put(28,8){\vector(1,1){4}}
\put(-12,12){\makebox(0,0)[cc]{$y_3$}} \put(-12,8){\vector(1,1){4}}
\put(88,32){\makebox(0,0)[cc]{$y_3$}} \put(88,28){\vector(1,1){4}}
\put(48,32){\makebox(0,0)[cc]{$y_1$}} \put(48,28){\vector(1,1){4}}
\put(8,32){\makebox(0,0)[cc]{$y_2$}} \put(8,28){\vector(1,1){4}}
\put(112,32){\makebox(0,0)[cc]{$z_1$}}
\put(112,28){\vector(-1,1){4}} \put(72,32){\makebox(0,0)[cc]{$z_2$}}
\put(72,28){\vector(-1,1){4}} \put(32,32){\makebox(0,0)[cc]{$z_3$}}
\put(32,28){\vector(-1,1){4}} \put(92,12){\makebox(0,0)[cc]{$z_1$}}
\put(92,8){\vector(-1,1){4}} \put(52,12){\makebox(0,0)[cc]{$z_2$}}
\put(52,8){\vector(-1,1){4}} \put(12,12){\makebox(0,0)[cc]{$z_3$}}
\put(12,8){\vector(-1,1){4}} \put(103,20){\makebox(0,0)[cc]{$z_1$}}
\put(100,22){\vector(0,-1){4}} \put(63,20){\makebox(0,0)[cc]{$z_2$}}
\put(60,22){\vector(0,-1){4}} \put(23,20){\makebox(0,0)[cc]{$z_3$}}
\put(20,22){\vector(0,-1){4}}
\end{picture}
\end{center} \caption{Identification $\varphi_3$ of faces of $\mathcal P_3$.}
\label{fig:m24}
\end{figure}
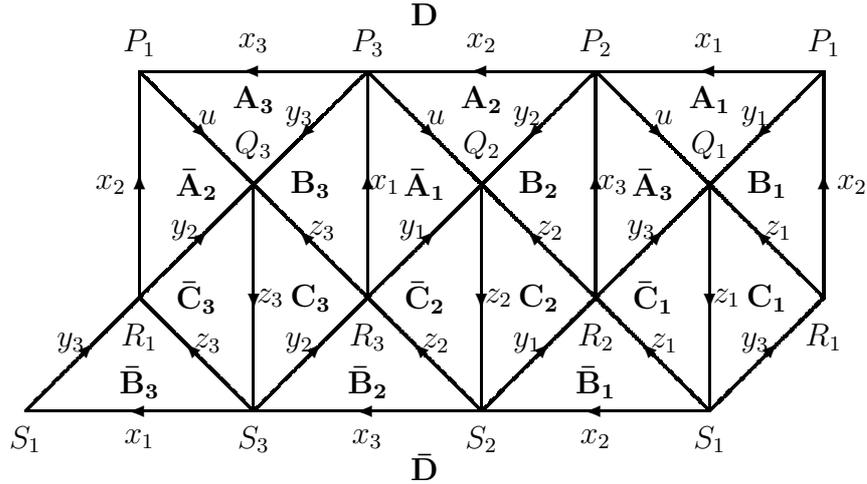

Denote by  $\mathcal P_3$ an icosahedron with all dihedral angles $2 \pi /3$. It is well-known that $\mathcal P_3$ can be realized in the  hyperbolic space $\mathbb H^3$.  Recall that $\mathcal P_3$ has 12 vertices, 30 edges and 20 faces. We will present $\mathcal P_3$ as in Fig.~\ref{fig:m24} where left and right sides, both denoted by $P_1 R_1 S_1$, are supposed to be identified.
The following pairwise identification $\varphi_3$ of the faces of $\mathcal P_3$, with the given ordering of vertices on the faces, can be found in  \cite{Everitt}:
\begin{equation}
\begin{array}{llll}
a_i  : {\bf A_i}  \to {\bf \bar{A}_i}  &  [P_i P_{i+1} Q_i  \to R_{i+2} P_{i+2} Q_{i+1}] ,  \quad &
b_i : {\bf B_i}  \to {\bf \bar{B}_i} & [R_i P_i Q_i  \to S_i S_{i+1} R_{i+1}] , \cr
c_i : {\bf C_i} \to {\bf \bar{C}_i} &  [S_i R_i Q_i \to R_{i+1} Q_i S_i], &
d : {\bf D} \to {\bf \bar{D}} &  [P_1 P_2 P_3 \to S_3 S_1 S_2] ,
\end{array}
\label{def:m24}
\end{equation}
where $i=1,2,3$ and all indices are taken mod $3$. Obviously, the face identification $\varphi_{3} = \{ a_{i}, b_{i}, c_{i}, d \}$ induces equivalent relations on the sets of vertices, edges, and faces of $\mathcal P_{3}$.  Since $\varphi_{3}$ can be realized as isometries of $\mathbb H^{3}$,  the quotient space $\mathcal P_3 / \varphi_3$ is a compact orientable hyperbolic 3-manifold, that was  denoted by $M_{24}$ in \cite{Cavicchioli, Cavicchioli2, Everitt}. It was shown in \cite{Cavicchioli2} that $M_{24}$ has the following interesting property: it is a 3-fold cyclic branched covering of the lens space  $L_{3,1}$ branched over a 2-component link.

\begin{figure}[h]
\begin{center}
\unitlength=0.7mm
\begin{picture}(200,80)(0,-5)
\thicklines
\qbezier(0,0)(0,0)(80,0) \qbezier(20,60)(20,60)(100,60)
\qbezier(100,0)(100,0)(180,0) \qbezier(120,60)(120,60)(200,60)
\qbezier(40,0)(40,0)(40,40) \qbezier(80,0)(80,0)(80,40)
\qbezier(140,0)(140,0)(140,40)  \qbezier(180,0)(180,0)(180,40)
\qbezier(20,20)(20,20)(20,60) \qbezier(60,20)(60,20)(60,60)
\qbezier(100,20)(100,20)(100,60) \qbezier(120,20)(120,20)(120,60)
\qbezier(160,20)(160,20)(160,60)  \qbezier(200,20)(200,20)(200,60)
\qbezier(40,40)(40,40)(20,20) \qbezier(40,40)(40,40)(20,60)
\qbezier(40,40)(40,40)(60,20)  \qbezier(40,40)(40,40)(60,60)
\qbezier(80,40)(80,40)(60,20) \qbezier(80,40)(80,40)(60,60)
\qbezier(80,40)(80,40)(100,20)  \qbezier(80,40)(80,40)(100,60)
\qbezier(140,40)(140,40)(120,20) \qbezier(140,40)(140,40)(120,60)
\qbezier(140,40)(140,40)(160,20)  \qbezier(140,40)(140,40)(160,60)
\qbezier(180,40)(180,40)(160,20) \qbezier(180,40)(180,40)(160,60)
\qbezier(180,40)(180,40)(200,20)  \qbezier(180,40)(180,40)(200,60)
\qbezier(0,0)(0,0)(20,20) \qbezier(20,20)(20,20)(40,0)
\qbezier(40,0)(40,0)(60,20)  \qbezier(60,20)(60,20)(80,0)
\qbezier(80,0)(80,0)(100,20) \qbezier(120,20)(120,20)(140,0)
\qbezier(160,20)(160,20)(140,0) \qbezier(160,20)(160,20)(180,0)
\qbezier(200,20)(200,20)(180,0) \qbezier(100,0)(100,0)(120,20)
\put(200,65){\makebox(0,0)[cc]{$P_1$}}
\put(160,65){\makebox(0,0)[cc]{$P_2$}}
\put(120,65){\makebox(0,0)[cc]{$P_3$}}
\put(100,65){\makebox(0,0)[cc]{$P_{n-1}$}}
\put(60,65){\makebox(0,0)[cc]{$P_n$}}
\put(20,65){\makebox(0,0)[cc]{$P_1$}}
\put(180,48){\makebox(0,0)[cc]{$Q_1$}}
\put(140,48){\makebox(0,0)[cc]{$Q_2$}}
\put(80,48){\makebox(0,0)[cc]{$Q_{n-1}$}}
\put(40,48){\makebox(0,0)[cc]{$Q_n$}}
\put(200,12){\makebox(0,0)[cc]{$R_1$}}
\put(160,12){\makebox(0,0)[cc]{$R_2$}}
\put(120,12){\makebox(0,0)[cc]{$R_3$}}
\put(100,12){\makebox(0,0)[cc]{$R_{n-1}$}}
\put(60,12){\makebox(0,0)[cc]{$R_n$}}
\put(20,12){\makebox(0,0)[cc]{$R_1$}}
\put(180,-5){\makebox(0,0)[cc]{$S_1$}}
\put(140,-5){\makebox(0,0)[cc]{$S_2$}}
\put(100,-5){\makebox(0,0)[cc]{$S_3$}}
\put(80,-5){\makebox(0,0)[cc]{$S_{n-1}$}}
\put(40,-5){\makebox(0,0)[cc]{$S_n$}}
\put(0,-5){\makebox(0,0)[cc]{$S_1$}}
\put(180,55){\makebox(0,0)[cc]{$\bf A_1$}}
\put(140,55){\makebox(0,0)[cc]{$\bf A_2$}}
\put(80,55){\makebox(0,0)[cc]{$\bf A_{n-1}$}}
\put(40,55){\makebox(0,0)[cc]{$\bf A_n$}}
\put(190,40){\makebox(0,0)[cc]{$\bf B_1$}}
\put(170,40){\makebox(0,0)[cc]{$\bf \bar{A}_n$}}
\put(150,40){\makebox(0,0)[cc]{$\bf B_2$}}
\put(130,40){\makebox(0,0)[cc]{$\bf \bar{A}_1$}}
\put(90,40){\makebox(0,0)[cc]{$\bf B_{n-1}$}}
\put(70,40){\makebox(0,0)[cc]{$\bf \bar{A}_{n-2}$}}
\put(50,40){\makebox(0,0)[cc]{$\bf B_n$}}
\put(30,40){\makebox(0,0)[cc]{$\bf \bar{A}_{n-1}$}}
\put(190,20){\makebox(0,0)[cc]{$\bf C_1$}}
\put(170,20){\makebox(0,0)[cc]{$\bf \bar{C}_1$}}
\put(150,20){\makebox(0,0)[cc]{$\bf C_2$}}
\put(130,20){\makebox(0,0)[cc]{$\bf \bar{C}_2$}}
\put(90,20){\makebox(0,0)[cc]{$\bf C_{n-1}$}}
\put(70,20){\makebox(0,0)[cc]{$\bf \bar{C}_{n-1}$}}
\put(50,20){\makebox(0,0)[cc]{$\bf C_n$}}
\put(30,20){\makebox(0,0)[cc]{$\bf \bar{C}_n$}}
\put(160,5){\makebox(0,0)[cc]{$\bf \bar{B}_1$}}
\put(120,5){\makebox(0,0)[cc]{$\bf \bar{B}_2$}}
\put(60,5){\makebox(0,0)[cc]{$\bf \bar{B}_{n-1}$}}
\put(20,5){\makebox(0,0)[cc]{$\bf \bar{B}_n$}}
\put(110,70){\makebox(0,0)[cc]{$\bf D$}}
\put(90,-10){\makebox(0,0)[cc]{$\bf \bar{D}$}}
\put(110,60){\makebox(0,0)[cc]{$\cdots$}}
\put(90,0){\makebox(0,0)[cc]{$\cdots$}}
\end{picture}
\end{center} \caption{Identification $\varphi_n$ of faces of $\mathcal P_n$.}
\label{fig:m24n}
\end{figure}
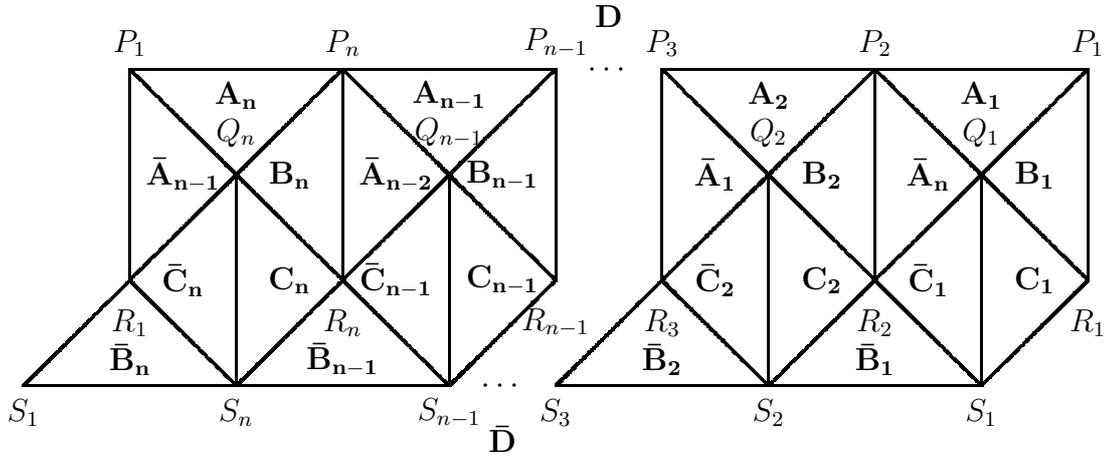

To generalize the construction of $M_{24}$, let us consider the complex $\mathcal P_n$, $n \geqslant 1$, having $4n$ vertices, $10n$ edges, and  $6n+2$ faces, presented in Fig.~\ref{fig:m24n}. In particular, $\mathcal P_3$ is the icosahedron as above. Define the pairwise identification $\varphi_n$ of faces of  $\mathcal P_n$ by formulae~(\ref{def:m24}) for $i=1, \ldots, n$, with  the following correction:
\begin{equation}
d : {\bf D} \to {\bf \bar{D}} \qquad  [P_1 P_2 \ldots P_{n-1} P_n \to S_3 S_4 \ldots S_1 S_2] ,
\end{equation}
and denote the corresponding  quotient space by~$M_{24}(n)$.
The following result was stated in \cite{Cavicchioli2} and~\cite{Kozlovskaya1, Kozlovskaya2}.
\begin{prop} \label{prop:m24}
For each $n \geqslant 1$ the quotient space $M_{24}(n) = \mathcal P_n / \varphi_n$ is a manifold.
\end{prop}

\dimo
Denote by $\sigma_k$ for $k=0,1,2,3$ the number of $k$-dimensional cells in $M_{24}(n)$. Obviously, $\sigma_3 = 1$. Moreover, $\sigma_2 = 3n+1$, since there are the following classes of equivalent faces: $\bf A_i \equiv {\bf \bar{A}_i}$, $\bf B_i \equiv {\bf \bar{B}_i}$, $\bf C_i \equiv {\bf \bar{C}_i}$, where  $i=1,\ldots,n$, and $\bf D \equiv {\bf \bar{D}}$.  Also, $\sigma_1 = 3n+1$, since all $1$-cells are separated in four types of equivalence classes:
\begin{eqnarray*}
({\rm I}_i) &  P_i P_{i+1} \xrightarrow{a_i} R_{i+2} P_{i+2} \xrightarrow{b_{i+2}} S_{i+2} S_{i+3} \xrightarrow{d^{-1}}  P_i P_{i+1} ;& \cr
({\rm II}_i) &  P_i Q_i \xrightarrow{a_i} R_{i+2} Q_{i+1} \xrightarrow{c_{i+1}^{-1}} S_{i+1} R_{i+1} \xrightarrow{b_{i}^{-1}} P_i Q_i ; & \cr
({\rm III}_i)  & Q_i R_i \xrightarrow{c_i} S_i Q_i \xrightarrow{c_i} R_{i+1} S_i \xrightarrow{b_i^{-1}} Q_i R_i   ; & \cr
({\rm IV}) & P_2 Q_1 \xrightarrow{a_1} P_3 Q_2 \xrightarrow{a_2} \ldots P_n Q_{n-1} \xrightarrow{a_{n-1}} P_1 Q_n \xrightarrow{a_n}  P_2 Q_1 , &
\end{eqnarray*}
where $i=1,\ldots,n$. It is easy to check that, by the action of $\varphi_n$, all vertices of $\mathcal P_n$ are equivalent, and so, $\sigma_0 = 1$. Thus, the Euler characteristic of the quotient space  is $\chi(M_{24}(n)) = 0$, and by \cite{Seifert-Threlfall}  $M_{24}(n)$ is a manifold.
\qed

We recall that a presentation for the fundamental group of a closed 3-manifold is geometric if it corresponds to a Heegaard diagram. The following presentation for the fundamental group of $M_{24}(n)$ was found in \cite{Kozlovskaya2}.

\begin{prop} \label{M24n-presentation}
The fundamental group of $M_{24}(n)$, $n\geqslant 1$, has the following geometric presentation:
\begin{equation}
\begin{gathered}
\pi_1 (M_{24}(n)) = \langle a_1, \ldots, a_n;  b_1, \ldots, b_n; c_1, \ldots, c_n;  d \quad | \quad  a_1 a_2 \ldots a_n = 1, \\
a_i b_{i+2} d^{-1} = 1, \qquad  a_i c_{i+1}^{-1} b_i^{-1} = 1, \qquad c_i^2 b_i^{-1}  = 1, \qquad i=1, \ldots, n \rangle.
\end{gathered} \label{eq2}
\end{equation}
\end{prop}

\dimo
An open Heegaard diagram for $M_{24}(n)$  arises from Fig.~\ref{fig:m24n} with discs corresponding to faces of $\mathcal P_n$ and segments of curves being dual to edges of $\mathcal P_n$. Thus, $\pi_1 (M_{24}(n))$ is generated by $a_1, \ldots, a_n$, $b_1,  \ldots, b_n$, $c_1, \ldots, c_n$, $d$. Obviously, passing along curves on the open Heegaard diagram we get defining relations as $({\rm I}_i)$, $({\rm II}_i)$, $({\rm III}_i)$ and $({\rm IV})$ in Proposition~\ref{prop:m24}: $a_i b_{i+2} d^{-1} = 1$, $a_i c_{i+1}^{-1} b_i^{-1} = 1$,  $c_i c_i b_i^{-1}=1$, and $a_1 a_2 \cdots a_n = 1$, respectively.
\qed

\begin{cor} \label{cor-m24n-presentation}
The fundamental group of $M_{24}(n)$, $n\geqslant 1$, has the following presentation:
\begin{equation}
\langle c_1, \ldots, c_n \quad | \quad  \prod_{j=1}^n c_j^3 = 1, \quad c_i^2 c_{i+1} c_{i+2}^2  = c_{i+1}^2 c_{i+2} c_{i+3}^2, \quad i=1, \ldots, n \rangle .
\end{equation}
\end{cor}

\dimo
Let us start with the group presentation given in Proposition~\ref{M24n-presentation}. Since $b_i = c_i^2$, we get $a_i = c_i^2 c_{i+1}$ for $i=1, \ldots n$, with indices taken mod $n$. Substituting these expressions into $a_1 \ldots a_n = 1$, we get $\prod_{j=1}^n (c_j^2 c_{j+1}) = 1$, and eliminating $d$ from $a_i b_{i+2} d^{-1} =1$, we get $d = c_i^2 c_{i+1} c_{i+2}^2$ for $i=1, \ldots n$, with indices taken mod $n$.
\qed

A presentation for $\pi_1(M_{24}(n))$, which is dual to (\ref{eq2}), was given in \cite{Cavicchioli2}:
\begin{equation}
\begin{gathered}
\pi_1 (M_{24}(n)) = \langle x_1, \ldots, x_n;  y_1, \ldots, y_n; z_1, \ldots, z_n;  u \quad | \quad  x_1 x_2 \ldots x_n = 1, \\
x_i u = y_i, \qquad x_i y_{i+2} = z_{i+2}, \qquad z_i^2 y_{i-1} = 1 , \qquad i=1, \ldots, n \rangle .
\end{gathered} \label{eq4}
\end{equation}
Here the relations correspond to the boundaries of the 2-faces, as shown for $M_{24}(3)$ in  Fig.~\ref{fig:m24}

\section{Covering properties and volumes of manifolds $M_{24}(n)$}
\label{sec2}

Let $M$, $M'$ be compact connected  orientable 3-manifolds and $L'$ a disjoint union of closed curves properly embedded in $M'$. Let $p : M \to M'$ denote a cyclic covering of $M'$ by $M$ branched over $L'$. We say that the branched covering $p$ is \emph{strongly-cyclic} if the stabilizer of each point of the singular set $p^{-1} (L')$ is the whole group of covering transformations.

The following property of $M_{24}(n)$ was observed in \cite{Cavicchioli2}. We will give a new proof, which, unlike the one in \cite{Cavicchioli2}, does not use results of \cite{Osborn-Stevens} and \cite{Stevens}.

\begin{thm} \label{theorem:m24n-covering}
For each $n \geqslant 2$, the manifold $M_{24}(n)$ is an $n$-fold strongly-cyclic branched covering of the lens space $L_{3,1}$, branched over a 2-component link. Moreover, $M_{24}(1)$ is the lens space~$L_{3,1}$.
\end{thm}

\dimo
Denote by $\rho_n$ the rotational symmetry of $\mathcal P_n$ sending $X_i$ to  $X_{i+1}$$, i=1, \ldots, n$, with indices taken mod $n$,  where $X$ belongs to the set of letters used for the notations of the vertices: $\{ P, Q, R, S\}$. This symmetry induces a cyclic symmetry of the quotient space $M_{24}(n) = \mathcal P_n / \varphi_n$, and we denote it by $\rho_n$, too. The quotient space $M_{24}(n) / \rho_n$ is an orbifold whose underlying manifold is $M_{24}(1)$.  Its singular set $\mathcal L$ consists of two components: $\mathcal L = \ell_1 \cup \ell_2$. According to the description of the equivalence classes of the edges given in Proposition~\ref{prop:m24}, the first component $\ell_1$ corresponds to the class ${\rm (IV)}$ of edges. The second component $\ell_2$ corresponds to the axis of rotation $\rho_n$.  Both components have singularity index $n$. Through the Heegaard diagram of $M_{24}(n) / \rho_n$, we understand the Heegaard diagram of $M_{24}(1)$ with information about the singular set $\mathcal L$ presented.

\begin{figure}[h]
\centering{
\includegraphics[height=6cm]{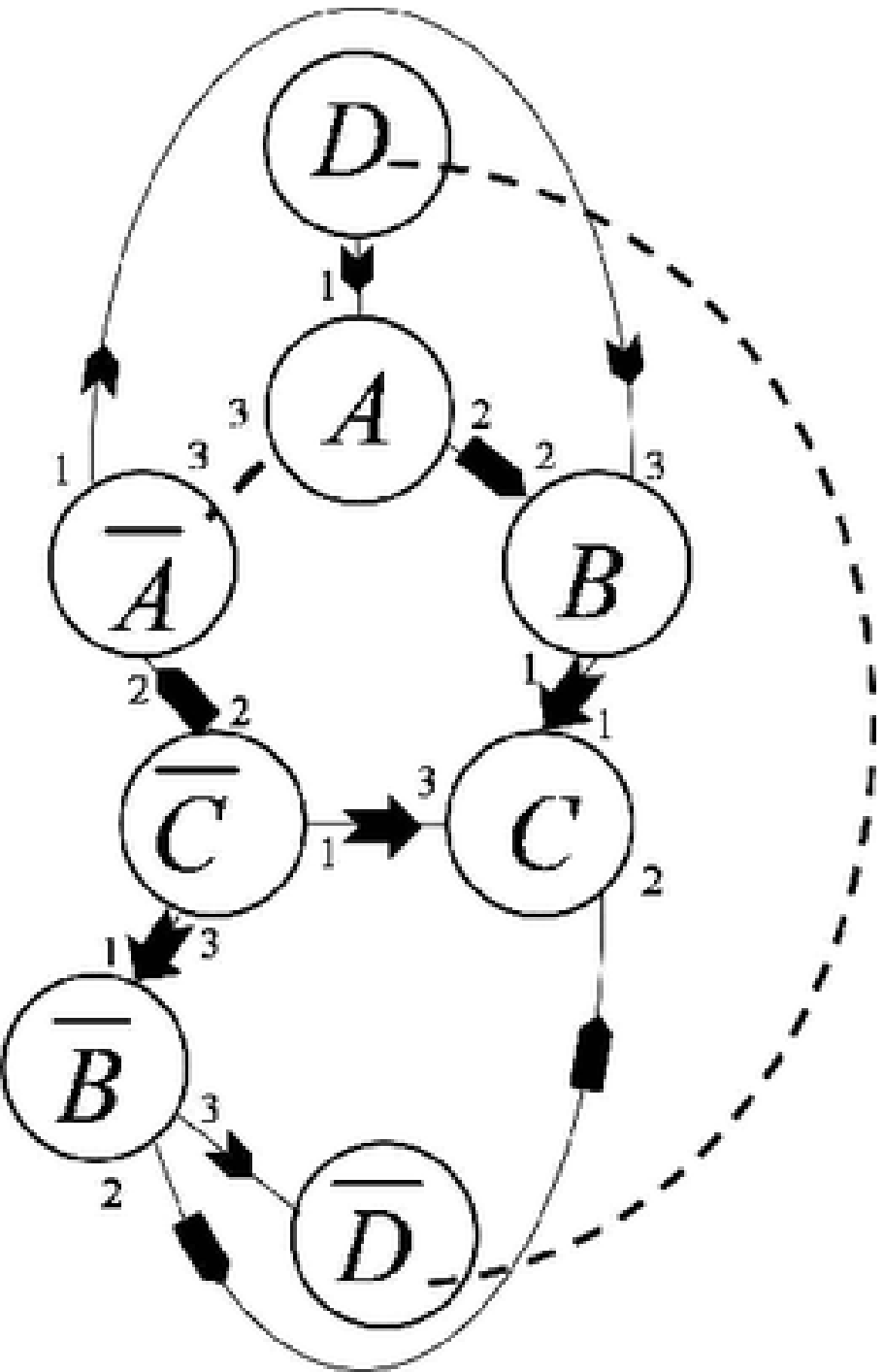}  \qquad
\includegraphics[height=6cm]{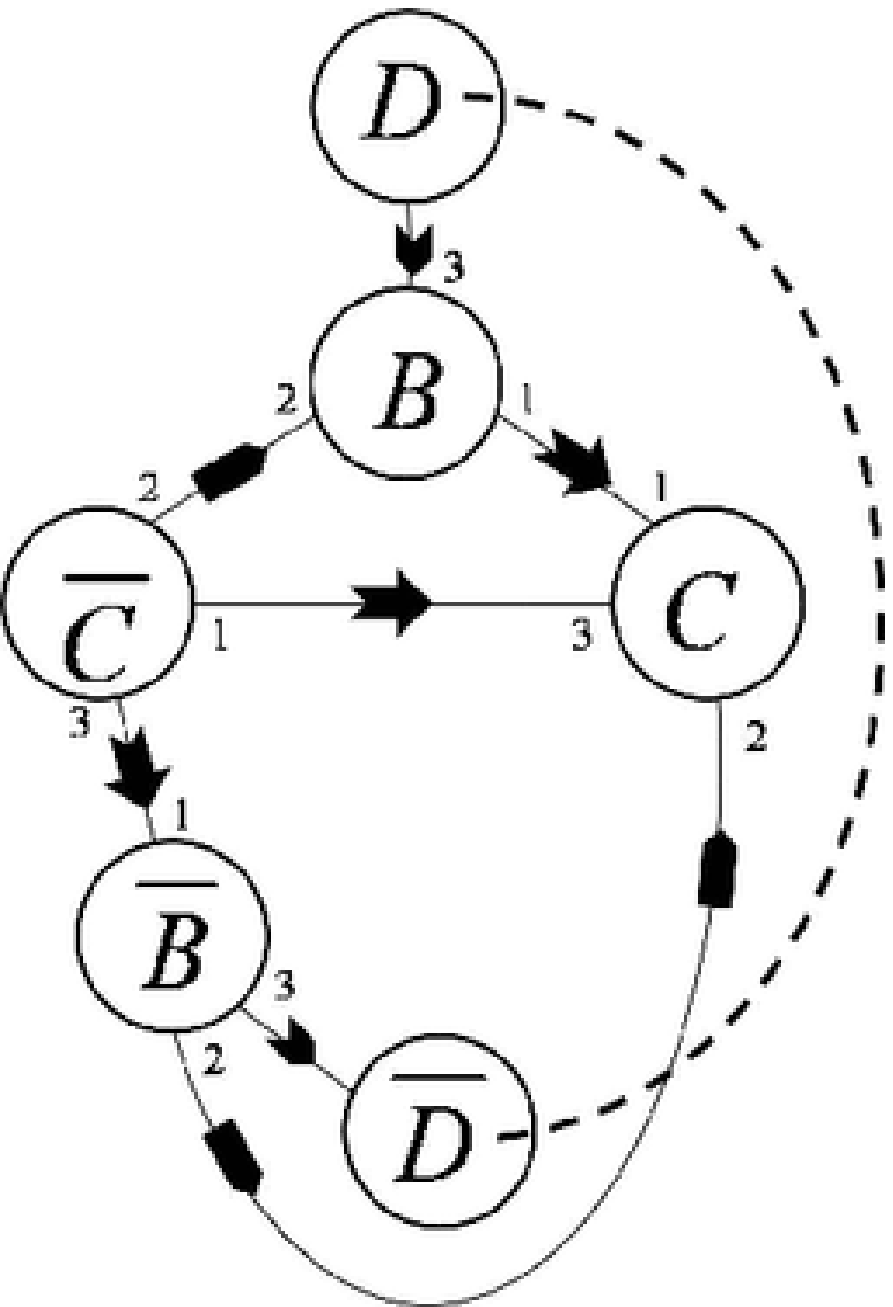}  \qquad
\includegraphics[height=4cm]{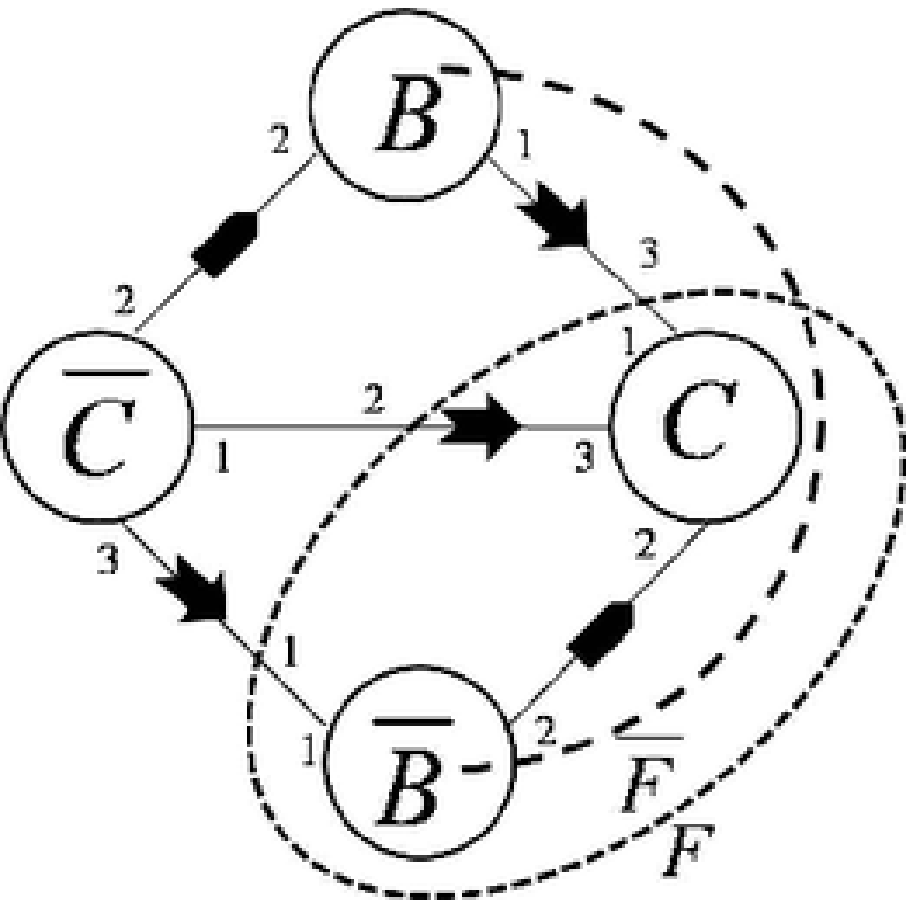}   \\
\includegraphics[height=4cm]{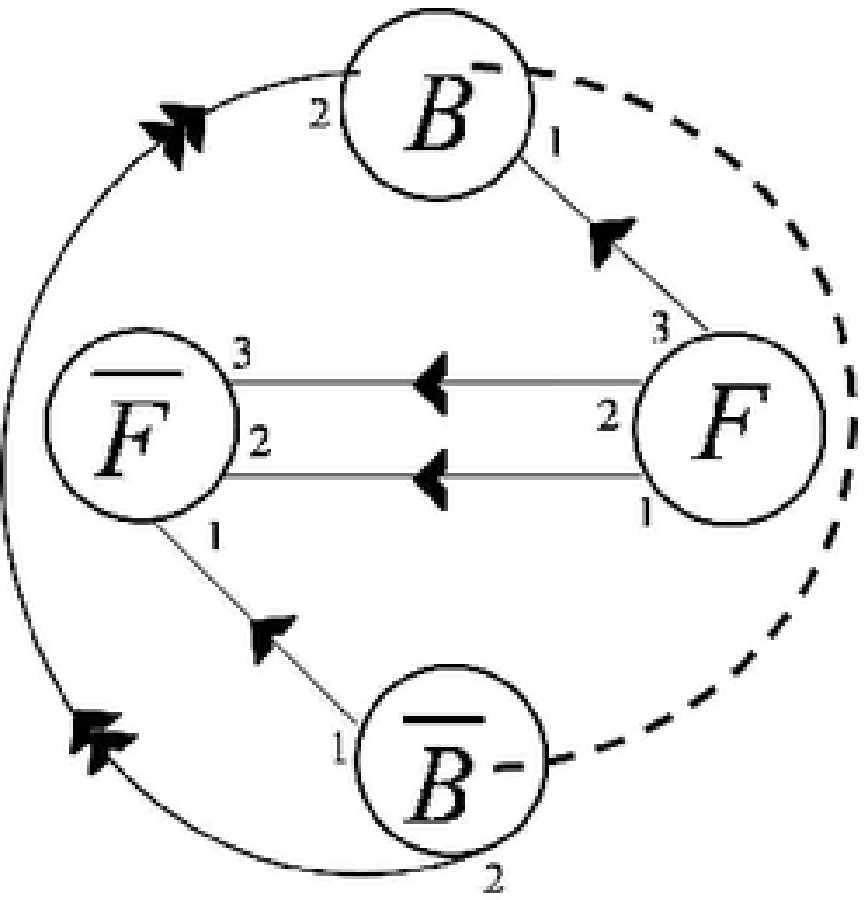}  \qquad  \qquad
\includegraphics[height=1.5cm]{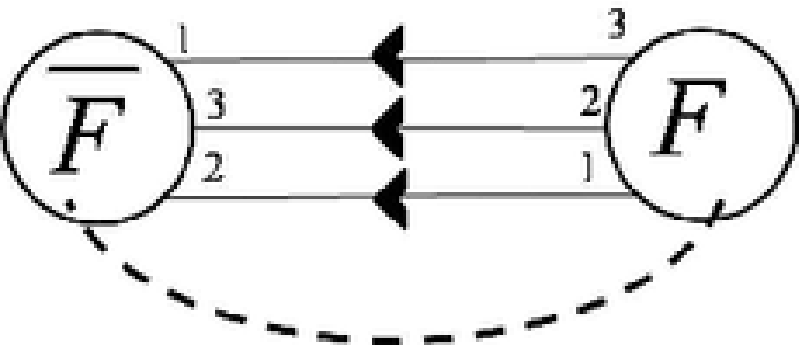}
\caption{Heegaard diagrams from $M_{24}(n) / \rho_n$ to $L_{3,1}$.} \label{fig:lens_from_m24}
}
\end{figure}

The equivalence transformations of Heegaard diagrams from $M_{24}(n) / \rho_n$ to $L_{3,1}$ are drawn in Fig.~\ref{fig:lens_from_m24}. Here the first component $\ell_{1}$ of the singular set is represented by a dashed segment connecting the discs $A$ and $\bar{A}$; the second component $\ell_{2}$ by a dashed segment, connecting the  discs $D$ and $\bar{D}$.  Both components have branching index $n$. Thus, the branched covering is strongly-cyclic.
At the first step we identify the discs  $A$ and $\bar{A}$, forming an 1-handle. The dashed segment (which also corresponds to $\ell_{1}$), connecting these discs, will give a 2-handle to glue up this 1-handle (i.e. they form a pair of complementary handles).  Thus, $\ell_{1}$ is a trivial knot; we are not drawing it in the next figures.  At the second step we  cancel the discs $D$ and $\bar{D}$, since they are connected only with the discs $B$ and $\bar{B}$, respectively, and the connecting segments are glued together. At the third step we cut along the curve shown by the dotted line to form a new pair of discs $F$ and $\bar{F}$ and then identify the discs $C$ and $\bar{C}$. After that we easily get a genus one Heegaard diagram, with the discs $F$ and $\bar{F}$, which is the standard diagram for the lens space $L_{3,1}$, where the dotted line represents~$\ell_{2}$.
\qed

\begin{prop}
$M_{24}(2)$ is the Seifert manifold  $(S^2; (3,1), (3,2), (3,2), (1,-1))$.
\end{prop}

\dimo
By Corollary~\ref{cor-m24n-presentation}, $\pi_1(M_{24}(2))= \langle c_1,c_2 \, | \,   c_1^3c_2^3 = 1, \,  c_1^2 c_2 c_1^2  = c_2^2 c_1 c_2^2 \rangle $. The second relation is equivalent to $c_2^{-2} c_1^2 c_2 c_1^2  = c_1 c_2^2$; after multiplication by $c_1$, we obtain $c_1^2 c_2^2 = c_1 c_2^{-2} c_1^2 c_2 c_1^2$ and by using the first relation, i.e. $c_1^3 = c_2^{-3}$, we have $c_1^2 c_2^2 = c_1^4 (c_2 c_1^2)^2$ and thus $c_1^{-2} c_2^2 = (c_2 c_1^2)^2$. Let us set $c_1=a,\ b=c_2c_1^2$ and so $c_2=ba^{-2}$. Then, we have the following presentation:
\begin{eqnarray*}
\pi_1(M_{24}(2)) & = & \langle a,b \quad | \quad  a^3(ba^{-2})^3 = 1, \quad a^{-2} (ba^{-2})^2  = b ^2 \rangle \cr
& = & \langle a,b \quad | \quad  a^3 a^2 b^2 ba^{-2} = 1, \quad (ba^{-2})^2  = a ^2 b ^2 \rangle \cr
& = & \langle a,b \quad | \quad  a^3 b^3 = 1, \quad (ba^{-2})^2  = a ^2 b ^2 \rangle .
\end{eqnarray*}

Consider the Seifert manifold $M = (S^2; (3,1), (3,2), (3,2), (1,-1)) = (S^2; (3,1), (3,2), (3,-1))$.
The standard presentation of $\pi_1(M)$ is as following (see \cite{Orlik}):
\begin{equation*}
\langle x, y, z, h \quad | \quad  x y z = 1,
\quad x h  = h x, \quad y h = h y, \quad z h = hz, \quad x^3 h =1, \quad y^3 h^2 =1, \quad z^3 h^{-1} = 1 \rangle .
\end{equation*}
We obtain $h = z^3$ from the last relation, $y = x^{-1} z^{-1}$ from the first one and we get
\begin{equation*}
\pi_1(M) = \langle x, z \quad | \quad  x z^3 = z^3 x, \quad x^{-1} z^3 = z^3 x^{-1},  \quad x^3 z^3  = 1, \quad z^6 = (x^{-1} z^{-1})^{-3} \rangle .
\end{equation*}
Since the first two relations come from the third, we have
\begin{eqnarray*}
\pi_1(M)  = &   \langle x, z \quad | \quad   x^3 z^3  = 1, \quad z^6 = (zx)^3 \rangle  
& =  \quad \langle x, z \quad | \quad   x^3 z^3  = 1, \quad z^5 = x z x z x \rangle \cr
= &    \langle x, z \quad | \quad   x^3 z^3  = 1, \quad z^2 = x z x z x^4 \rangle
& =  \quad \langle x, z \quad | \quad   x^3 z^3  = 1, \quad z^2 x^2 = x z x z^{-5} \rangle \cr
= &   \langle x, z \quad | \quad   x^3 z^3  = 1, \quad z^2 x^2 = x z^{-2} x z^{-2} \rangle
& = \quad \langle x, z \quad | \quad   x^3 z^3  = 1, \quad z^2 x^2 = (x z^{-2})^2 \rangle ,
\end{eqnarray*}
which is the same presentation as above, i.e. the fundamental groups of $M_{24}(2)$ and $M$ are isomorphic.

As a consequence, note that, since $M$ is irreducible, $M_{24}(2)$ is irreducible, too. Furthermore, since $M$ is a large Seifert manifold, $\pi_1(M_{24}(2))=\pi_1(M)$ contains an infinite cyclic normal subgroup (generated by $h$, see \cite{Orlik}). Therefore, by a result of \cite{Casson-Jungreis} and \cite{Gabai} (see also \cite{Scott}), $M_{24}(2)$ is Seifert fibered and, consequently, $M_{24}(2)$ and $M$ are homeomorphic.
\qed

Since $M_{24}(2)$ is $(S^2; (3,1), (3,2), (3,2), (1,-1))$, it admits Nil geometry and its first homology group is $\mathbb Z_3 \oplus \mathbb Z_6$. Moreover, $M_{24}(2)$ can be obtained by Dehn surgeries with parameters $(-1,2)$, $(-2,1)$, $(-7,2)$ on the link
$6^3_1$ (\emph{chain link}), as well as with parameters $(6,1)$, $(1,1)$, $(3,1)$ on $6^3_2$ (\emph{Borromean rings}).

Theorem~3.1 from \cite{Cavicchioli2} states that for $n \geqslant 3$ the manifolds $M_{24}(n)$ are hyperbolic, although for $n > 3$ the authors do not present explicitly any proof of hyperbolicity. Moreover, the authors give the following volume formula:  $\textrm{\rm vol } M_{24} (n) = (n/3) \cdot (4.686034274\ldots)$. However, even for small $n > 3$ the above formula turns out to be wrong. It would be right if the polyhedron $\mathcal P_n$ could be obtained by gluing isometrically  $n$ copies of the $1/3$--piece of the hyperbolic $2\pi / 3$--icosahedron $\mathcal P_3$. This is obviously not true, since the dihedral angle  around the image of the axis of rotation $\rho_n$ in $\mathcal P_n / \rho_n$ must be equal to $2 \pi / n$, that is, it differs from $2 \pi / 3$ if $n > 3$.
The correct values of $\textrm{\rm vol } M_{24}(n)$ are presented below.

\begin{prop} \label{prop1.4}
The hyperbolic volumes and the first homology groups of $M_{24}(n)$, for $3\leqslant n\leqslant 6$, are as
follows:
$$
\begin{tabular}{|c|c|c|}
\hline {\rm manifold} & {\rm volume} & {\rm homology group} \cr
\hline  $M_{24}(3)$ & {\rm 4.686034273803\ldots} & $\mathbb Z_9$ \cr
\hline $M_{24}(4)$ & {\rm 9.702341514665\ldots} & $\mathbb Z_3 \oplus \mathbb Z_{12}$ \cr
\hline $M_{24}(5)$ & {\rm 14.319926985892\ldots} & $\mathbb Z_5 \oplus \mathbb Z_5 \oplus \mathbb Z_{15}$ \cr
\hline $M_{24}(6)$ & {\rm 18.649157163789\ldots} & $\mathbb Z_3 \oplus \mathbb Z_9 \oplus \mathbb Z_{18}$ \cr \hline
\end{tabular}
$$
\end{prop}

\dimo
Results are obtained by using the computer program \emph{Recognizer} .
\qed

\section{The family of manifolds $M_{25}(n)$}
\label{sec3}

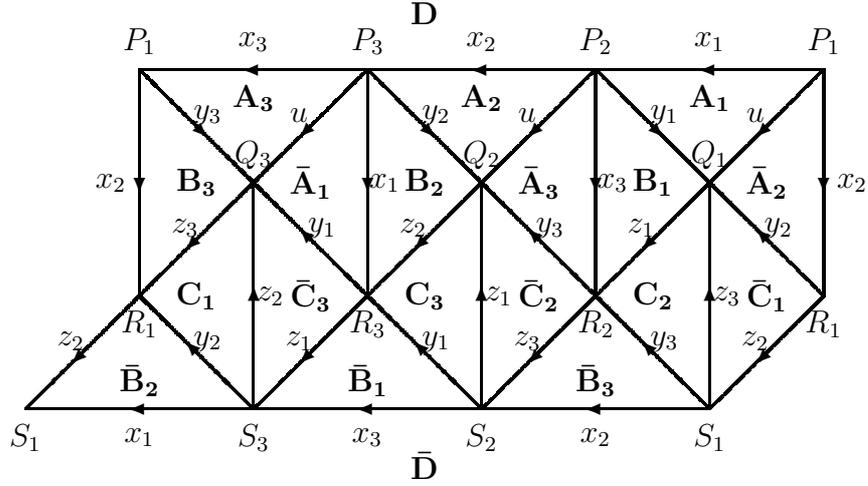
\begin{figure}[h]
\begin{center}
\unitlength=0.75mm
\begin{picture}(140,80)(-20,-5)
\thicklines
\qbezier(-20,0)(-20,0)(100,0) \qbezier(0,60)(0,60)(120,60)
\qbezier(20,0)(20,0)(20,40) \qbezier(60,0)(60,0)(60,40)
\qbezier(100,0)(100,0)(100,40)
\qbezier(0,20)(0,20)(0,60) \qbezier(40,20)(40,20)(40,60)
\qbezier(80,20)(80,20)(80,60) \qbezier(120,20)(120,20)(120,60)
\qbezier(20,40)(20,40)(0,20) \qbezier(20,40)(20,40)(0,60)
\qbezier(20,40)(20,40)(40,20)  \qbezier(20,40)(20,40)(40,60)
\qbezier(60,40)(60,40)(40,20) \qbezier(60,40)(60,40)(40,60)
\qbezier(60,40)(60,40)(80,20)  \qbezier(60,40)(60,40)(80,60)
\qbezier(100,40)(100,40)(80,20) \qbezier(100,40)(100,40)(80,60)
\qbezier(100,40)(100,40)(120,20)  \qbezier(100,40)(100,40)(120,60)
\qbezier(0,20)(0,20)(-20,0)
\qbezier(0,20)(0,20)(20,0) \qbezier(40,20)(40,20)(20,0)
\qbezier(40,20)(40,20)(60,0)  \qbezier(80,20)(80,20)(60,0)
\qbezier(80,20)(80,20)(100,0) \qbezier(120,20)(120,20)(100,0)
\put(120,65){\makebox(0,0)[cc]{$P_1$}}
\put(80,65){\makebox(0,0)[cc]{$P_2$}}
\put(40,65){\makebox(0,0)[cc]{$P_3$}}
\put(0,65){\makebox(0,0)[cc]{$P_1$}}
\put(100,45){\makebox(0,0)[cc]{$Q_1$}}
\put(60,45){\makebox(0,0)[cc]{$Q_2$}}
\put(20,45){\makebox(0,0)[cc]{$Q_3$}}
\put(120,15){\makebox(0,0)[cc]{$R_1$}}
\put(80,15){\makebox(0,0)[cc]{$R_2$}}
\put(40,15){\makebox(0,0)[cc]{$R_3$}}
\put(0,15){\makebox(0,0)[cc]{$R_1$}}
\put(100,-5){\makebox(0,0)[cc]{$S_1$}}
\put(60,-5){\makebox(0,0)[cc]{$S_2$}}
\put(20,-5){\makebox(0,0)[cc]{$S_3$}}
\put(-20,-5){\makebox(0,0)[cc]{$S_1$}}
\put(100,55){\makebox(0,0)[cc]{$\bf A_1$}}
\put(60,55){\makebox(0,0)[cc]{$\bf A_2$}}
\put(20,55){\makebox(0,0)[cc]{$\bf A_3$}}
\put(110,40){\makebox(0,0)[cc]{$\bf \bar{A}_2$}}
\put(90,40){\makebox(0,0)[cc]{$\bf B_1$}}
\put(70,40){\makebox(0,0)[cc]{$\bf \bar{A}_3$}}
\put(50,40){\makebox(0,0)[cc]{$\bf B_2$}}
\put(30,40){\makebox(0,0)[cc]{$\bf \bar{A}_1$}}
\put(10,40){\makebox(0,0)[cc]{$\bf B_3$}}
\put(110,20){\makebox(0,0)[cc]{$\bf \bar{C}_1$}}
\put(90,20){\makebox(0,0)[cc]{$\bf C_2$}}
\put(70,20){\makebox(0,0)[cc]{$\bf \bar{C}_2$}}
\put(50,20){\makebox(0,0)[cc]{$\bf C_3$}}
\put(30,20){\makebox(0,0)[cc]{$\bf \bar{C}_3$}}
\put(10,20){\makebox(0,0)[cc]{$\bf C_1$}}
\put(80,5){\makebox(0,0)[cc]{$\bf \bar{B}_3$}}
\put(40,5){\makebox(0,0)[cc]{$\bf \bar{B}_1$}}
\put(0,5){\makebox(0,0)[cc]{$\bf \bar{B}_2$}}
\put(50,70){\makebox(0,0)[cc]{$\bf D$}}
\put(50,-10){\makebox(0,0)[cc]{$\bf \bar{D}$}}
\put(100,65){\makebox(0,0)[cc]{$x_1$}}
\put(102,60){\vector(-1,0){4}}
\put(60,65){\makebox(0,0)[cc]{$x_2$}}
\put(62,60){\vector(-1,0){4}}
\put(20,65){\makebox(0,0)[cc]{$x_3$}}
\put(22,60){\vector(-1,0){4}}
\put(80,-5){\makebox(0,0)[cc]{$x_2$}}
\put(82,0){\vector(-1,0){4}}
\put(40,-5){\makebox(0,0)[cc]{$x_3$}}
\put(42,0){\vector(-1,0){4}}
\put(0,-5){\makebox(0,0)[cc]{$x_1$}}
\put(2,0){\vector(-1,0){4}}
\put(125,40){\makebox(0,0)[cc]{$x_2$}}
\put(120,42){\vector(0,-1){4}}
\put(83,40){\makebox(0,0)[cc]{$x_3$}}
\put(80,42){\vector(0,-1){4}}
\put(43,40){\makebox(0,0)[cc]{$x_1$}}
\put(40,42){\vector(0,-1){4}}
\put(-5,40){\makebox(0,0)[cc]{$x_2$}}
\put(0,42){\vector(0,-1){4}}
\put(92,52){\makebox(0,0)[cc]{$y_1$}}
\put(88,52){\vector(1,-1){4}}
\put(52,52){\makebox(0,0)[cc]{$y_2$}}
\put(48,52){\vector(1,-1){4}}
\put(12,52){\makebox(0,0)[cc]{$y_3$}}
\put(8,52){\vector(1,-1){4}}
\put(108,52){\makebox(0,0)[cc]{$u$}}
\put(112,52){\vector(-1,-1){4}}
\put(68,52){\makebox(0,0)[cc]{$u$}}
\put(72,52){\vector(-1,-1){4}}
\put(28,52){\makebox(0,0)[cc]{$u$}}
\put(32,52){\vector(-1,-1){4}}
\put(108,12){\makebox(0,0)[cc]{$z_2$}}
\put(112,12){\vector(-1,-1){4}}
\put(68,12){\makebox(0,0)[cc]{$z_3$}}
\put(72,12){\vector(-1,-1){4}}
\put(28,12){\makebox(0,0)[cc]{$z_1$}}
\put(32,12){\vector(-1,-1){4}}
\put(-12,12){\makebox(0,0)[cc]{$z_2$}}
\put(-8,12){\vector(-1,-1){4}}
\put(88,32){\makebox(0,0)[cc]{$z_1$}}
\put(92,32){\vector(-1,-1){4}}
\put(48,32){\makebox(0,0)[cc]{$z_2$}}
\put(52,32){\vector(-1,-1){4}}
\put(8,32){\makebox(0,0)[cc]{$z_3$}}
\put(12,32){\vector(-1,-1){4}}
\put(112,32){\makebox(0,0)[cc]{$y_2$}}
\put(112,28){\vector(-1,1){4}}
\put(72,32){\makebox(0,0)[cc]{$y_3$}}
\put(72,28){\vector(-1,1){4}}
\put(32,32){\makebox(0,0)[cc]{$y_1$}}
\put(32,28){\vector(-1,1){4}}
\put(92,12){\makebox(0,0)[cc]{$y_3$}}
\put(92,8){\vector(-1,1){4}}
\put(52,12){\makebox(0,0)[cc]{$y_1$}}
\put(52,8){\vector(-1,1){4}}
\put(12,12){\makebox(0,0)[cc]{$y_2$}}
\put(12,8){\vector(-1,1){4}}
\put(103,20){\makebox(0,0)[cc]{$z_3$}}
\put(100,18){\vector(0,1){4}}
\put(63,20){\makebox(0,0)[cc]{$z_1$}}
\put(60,18){\vector(0,1){4}}
\put(23,20){\makebox(0,0)[cc]{$z_2$}}
\put(20,18){\vector(0,1){4}}
\end{picture}
\end{center} \caption{Identification $\psi_3$ of faces of ${\mathcal P}_3$.} \label{fig:m25}
\end{figure}

The following pairwise identification $\psi_3$ of faces of $\mathcal P_3$, with notations according to Fig.~\ref{fig:m25},  can be found in  \cite{Everitt}:
\begin{equation}
\begin{array}{llll}
a_i : {\bf A_i}  \to {\bf \bar{A}_i} & [ P_i P_{i+1} Q_i  \to P_{i+2} R_{i+2} Q_{i+2}]  , &
b_i : {\bf B_i}  \to {\bf \bar{B}_i} & [ Q_i R_{i+1} P_{i+1} \to R_{i+2} S_{i+2} S_{i+1} ]  , \cr
c_i : {\bf C_i} \to {\bf \bar{C}_i} & [ Q_{i-1} R_i S_{i-1} \to S_i Q_i R_i ] , &
d : {\bf D} \to {\bf \bar{D}}  & [P_1 P_2 P_3 \to S_3 S_1 S_2] ,
\end{array}
\label{def:m25}
\end{equation}
where $i=1,2,3$ and all indices are taken mod $3$.
The quotient space $\mathcal P_3 / \psi_3$ is a compact orientable hyperbolic 3-manifold denoted by $M_{25}$ in \cite{Cavicchioli, Cavicchioli2, Everitt}. It was shown in \cite{Cavicchioli2} that $M_{25}$ is a 3-fold cyclic branched covering of the lens space $L_{3,1}$ branched over a 2-component link.

\begin{figure}[h]
\begin{center}
\unitlength=0.7mm
\begin{picture}(200,80)(0,-5)
\thicklines
\qbezier(0,0)(0,0)(80,0) \qbezier(20,60)(20,60)(100,60)
\qbezier(100,0)(100,0)(180,0) \qbezier(120,60)(120,60)(200,60)
\qbezier(40,0)(40,0)(40,40) \qbezier(80,0)(80,0)(80,40)
\qbezier(140,0)(140,0)(140,40)  \qbezier(180,0)(180,0)(180,40)
\qbezier(20,20)(20,20)(20,60) \qbezier(60,20)(60,20)(60,60)
\qbezier(100,20)(100,20)(100,60) \qbezier(120,20)(120,20)(120,60)
\qbezier(160,20)(160,20)(160,60)  \qbezier(200,20)(200,20)(200,60)
\qbezier(40,40)(40,40)(20,20) \qbezier(40,40)(40,40)(20,60)
\qbezier(40,40)(40,40)(60,20)  \qbezier(40,40)(40,40)(60,60)
\qbezier(80,40)(80,40)(60,20) \qbezier(80,40)(80,40)(60,60)
\qbezier(80,40)(80,40)(100,20)  \qbezier(80,40)(80,40)(100,60)
\qbezier(140,40)(140,40)(120,20) \qbezier(140,40)(140,40)(120,60)
\qbezier(140,40)(140,40)(160,20)  \qbezier(140,40)(140,40)(160,60)
\qbezier(180,40)(180,40)(160,20) \qbezier(180,40)(180,40)(160,60)
\qbezier(180,40)(180,40)(200,20)  \qbezier(180,40)(180,40)(200,60)
\qbezier(0,0)(0,0)(20,20) \qbezier(20,20)(20,20)(40,0)
\qbezier(40,0)(40,0)(60,20)  \qbezier(60,20)(60,20)(80,0)
\qbezier(80,0)(80,0)(100,20) \qbezier(120,20)(120,20)(140,0)
\qbezier(160,20)(160,20)(140,0) \qbezier(160,20)(160,20)(180,0)
\qbezier(200,20)(200,20)(180,0) \qbezier(100,0)(100,0)(120,20)
\put(200,65){\makebox(0,0)[cc]{$P_1$}}
\put(160,65){\makebox(0,0)[cc]{$P_2$}}
\put(120,65){\makebox(0,0)[cc]{$P_3$}}
\put(100,65){\makebox(0,0)[cc]{$P_{n-1}$}}
\put(60,65){\makebox(0,0)[cc]{$P_n$}}
\put(20,65){\makebox(0,0)[cc]{$P_1$}}
\put(180,48){\makebox(0,0)[cc]{$Q_1$}}
\put(140,48){\makebox(0,0)[cc]{$Q_2$}}
\put(80,48){\makebox(0,0)[cc]{$Q_{n-1}$}}
\put(40,48){\makebox(0,0)[cc]{$Q_n$}}
\put(200,12){\makebox(0,0)[cc]{$R_1$}}
\put(160,12){\makebox(0,0)[cc]{$R_2$}}
\put(120,12){\makebox(0,0)[cc]{$R_3$}}
\put(100,12){\makebox(0,0)[cc]{$R_{n-1}$}}
\put(60,12){\makebox(0,0)[cc]{$R_n$}}
\put(20,12){\makebox(0,0)[cc]{$R_1$}}
\put(180,-5){\makebox(0,0)[cc]{$S_1$}}
\put(140,-5){\makebox(0,0)[cc]{$S_2$}}
\put(100,-5){\makebox(0,0)[cc]{$S_3$}}
\put(80,-5){\makebox(0,0)[cc]{$S_{n-1}$}}
\put(40,-5){\makebox(0,0)[cc]{$S_n$}}
\put(0,-5){\makebox(0,0)[cc]{$S_1$}}
\put(180,55){\makebox(0,0)[cc]{$\bf A_1$}}
\put(140,55){\makebox(0,0)[cc]{$\bf A_2$}}
\put(80,55){\makebox(0,0)[cc]{$\bf A_{n-1}$}}
\put(40,55){\makebox(0,0)[cc]{$\bf A_n$}}
\put(190,40){\makebox(0,0)[cc]{$\bf \bar{A}_{n-1}$}}
\put(170,40){\makebox(0,0)[cc]{$\bf B_1$}}
\put(150,40){\makebox(0,0)[cc]{$\bf \bar{A}_n$}}
\put(130,40){\makebox(0,0)[cc]{$\bf B_2$}}
\put(90,40){\makebox(0,0)[cc]{$\bf \bar{A}_{n-3}$}}
\put(70,40){\makebox(0,0)[cc]{$\bf B_{n-1}$}}
\put(50,40){\makebox(0,0)[cc]{$\bf \bar{A}_{n-2}$}}
\put(30,40){\makebox(0,0)[cc]{$\bf B_n$}}
\put(190,20){\makebox(0,0)[cc]{$\bf \bar{C}_1$}}
\put(170,20){\makebox(0,0)[cc]{$\bf C_2$}}
\put(150,20){\makebox(0,0)[cc]{$\bf \bar{C}_2$}}
\put(130,20){\makebox(0,0)[cc]{$\bf C_3$}}
\put(90,20){\makebox(0,0)[cc]{$\bf \bar{C}_{n-1}$}}
\put(70,20){\makebox(0,0)[cc]{$\bf C_n$}}
\put(50,20){\makebox(0,0)[cc]{$\bf \bar{C}_n$}}
\put(30,20){\makebox(0,0)[cc]{$\bf C_1$}}
\put(160,5){\makebox(0,0)[cc]{$\bf \bar{B}_n$}}
\put(120,5){\makebox(0,0)[cc]{$\bf \bar{B}_1$}}
\put(60,5){\makebox(0,0)[cc]{$\bf \bar{B}_{n-2}$}}
\put(20,5){\makebox(0,0)[cc]{$\bf \bar{B}_{n-1}$}}
\put(110,70){\makebox(0,0)[cc]{$\bf D$}}
\put(90,-10){\makebox(0,0)[cc]{$\bf \bar{D}$}}
\put(110,60){\makebox(0,0)[cc]{$\cdots$}}
\put(90,0){\makebox(0,0)[cc]{$\cdots$}}
\end{picture}
\end{center} \caption{Identification $\psi_n$ of faces of ${\mathcal P}_n$.} \label{fig:m25n}
\end{figure}
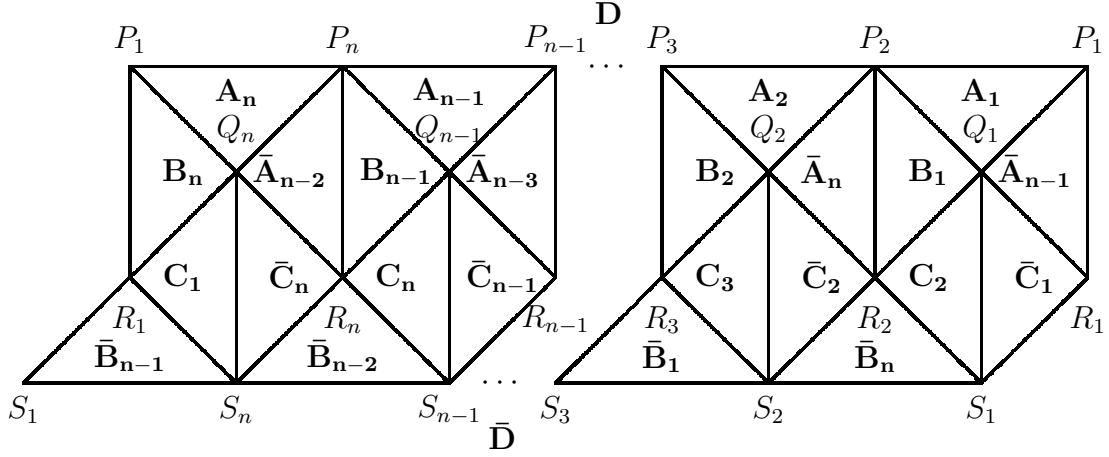

As one can see from Fig.~\ref{fig:m25}, the boundaries of the faces of $\mathcal P_{3}$ are in correspondence with the following relations:
\begin{equation}
x_{1} x_{2}  x_{3}  = 1, \quad x_i y_i = u, \quad y_i z_i  = x_{i-1}, \quad z_{i-1} z_i  = y_{i-1} , \qquad i=1, 2, 3 .
\end{equation}
To generalize the construction of $M_{25}$, it is natural to consider the complex $\mathcal P_n$, $n \geqslant 1$, pictured in Fig.~\ref{fig:m25n}, and define the pairwise identification $\psi_n$ of the faces of $\mathcal P_n$ by formulae~(\ref{def:m25}) for $i=1, \ldots, n$, with the following correction:
\begin{equation}
d : {\bf D} \to {\bf \bar{D}} \qquad  [P_1 P_2 \ldots P_{n-1} P_n \to S_3 S_4 \ldots S_1 S_2] .
\end{equation}
This is equivalent to the generalization from $M_{25}$ to $M_{25}(n)$ considered in \cite{Cavicchioli2}, were the face identifications of $\mathcal P_{n}$ are defined by the boundary relations corresponding to the defining relations for the following group presentation:
\begin{equation}
\begin{gathered}
G (n) = \langle x_1, \ldots, x_n;  y_1, \ldots, y_n; z_1, \ldots, z_n;  u \quad | \quad  x_1 x_2 \ldots x_n = 1, \\
x_i y_i = u, \qquad y_i z_i  = x_{i-1}, \qquad z_{i-1} z_i  = y_{i-1} , \qquad i=1, \ldots, n \rangle .
\end{gathered}   \label{group-G}
\end{equation}
Let us denote the corresponding  quotient space $\mathcal P_{n} /  \psi_{n}$ by $M_{25}(n)$.

\begin{prop} \label{prop:m25}
For each $n \geqslant 1$, the quotient space $M_{25} (n) = \mathcal P_n / \psi_n$ is a manifold.
\end{prop}

\dimo
Let  $\sigma_k$ be the number of $k$-cells in $M_{25}(n)$, $k=0,1,2,3$. Obviously,  $\sigma_3 = 1$. Moreover, $\sigma_2 = 3n+1$, since there are the following classes of equivalent faces: $\bf A_i \equiv {\bf \bar{A}_i}$, $\bf B_i \equiv {\bf \bar{B}_i}$, $\bf C_i \equiv {\bf \bar{C}_i}$, where  $i=1,\ldots,n$, and $\bf D \equiv {\bf \bar{D}}$. We have the following three types of classes of equivalent edges, with $n$ classes of each type:
\begin{eqnarray*}
({\rm I}_i) &  P_i P_{i+1} \xrightarrow{a_i} P_{i+2} R_{i+2} \xrightarrow{b_{i+1}} S_{i+2} S_{i+3} \xrightarrow{d^{-1}}  P_i P_{i+1}; & \cr ({\rm II}_i) &  P_{i+1} Q_i \xrightarrow{a_i} R_{i+2} Q_{i+2} \xrightarrow{c_{i+2}^{-1}} S_{i+1} R_{i+2} \xrightarrow{b_i^{-1}} P_{i+1} Q_i ; & \cr ({\rm III}_i)  & Q_{i-1}
R_i \xrightarrow{c_i} S_i Q_i  \xrightarrow{c_{i+1}} R_{i+1} S_{i+1}
\xrightarrow{b_{i-1}^{-1}} Q_{i-1} R_i  ; &
\end{eqnarray*}
where $i=1,\ldots,n$. Remark that $P_i Q_i \xrightarrow{a_i} P_{i+2} Q_{i+2}$. Thus, if $n$ is odd, then the set of edges $\{ P_1 Q_1 , \ldots P_n Q_n \}$ will form one class of equivalent edges:
$$
P_1 Q_1 \xrightarrow{a_1} P_{3} Q_{3} \xrightarrow{a_3} \ldots P_n Q_n \xrightarrow{a_n} P_2 Q_2 \xrightarrow{a_2} \ldots P_{n-1} Q_{n-1} \xrightarrow{a_{n-1}} P_1 Q_1  .
$$
Hence $\sigma_1 = 3n+1$. Moreover, in this case all vertices of $\mathcal P_n$ are equivalent, so $\sigma_0 = 1$. If  $n$ is even, we get two classes of equivalent edges:
$$
P_1 Q_1 \xrightarrow{a_1} P_{3} Q_{3} \xrightarrow{a_3} \ldots P_{n-1} Q_{n-1} \xrightarrow{a_{n-1}} P_1 Q_1; \qquad  P_2 Q_2 \xrightarrow{a_2}  P_4 Q_4 \xrightarrow{a_4} \ldots P_n Q_n \xrightarrow{a_n}  P_2 Q_2 .
$$
Hence $\sigma_1 = 3n+2$. Moreover, in this case  all vertices of $\mathcal P_n$ are separated in two classes of  equivalence, so $\sigma_0 = 2$. Thus, in both cases, the Euler characteristic of the quotient space is $\chi(M_{25}(n)) = 0$, and so $M_{25}(n)$ is a manifold.
\qed

\begin{prop}  \label{proposition-m25n-presentation}
The fundamental group of $M_{25}(n)$, $n\geqslant 1$, has the following geometric presentation:
\begin{equation}
\begin{gathered}
\pi_1 (M_{25}(n)) = \langle a_1, \ldots, a_n;  b_1, \ldots, b_n; c_1, \ldots, c_n;  d \quad | \quad a_1 a_3 \ldots a_n a_2 \ldots a_{n-1} =1, \\
a_i b_{i+1} d^{-1} = 1,  \, a_i c_{i+2}^{-1} b_i^{-1} = 1, \, c_i c_{i+1} b_{i-1}^{-1}= 1, \quad i=1, \ldots, n \rangle .
\end{gathered} \label{eq6}
\end{equation}
if $n$ is odd, and
\begin{equation}
\begin{gathered}
\pi_1 (M_{25}(n)) = \langle a_1, \ldots, a_n;  b_1, \ldots, b_n; c_1, \ldots, c_n;  d \quad | \quad a_1 a_3 \ldots a_{n-1} = 1,  \quad a_2 a_4 \ldots a_n =1 ,  \\
a_i b_{i+1} d^{-1} = 1,  \, a_i c_{i+2}^{-1} b_i^{-1} = 1, \, c_i c_{i+1} b_{i-1}^{-1}= 1, \quad i=1, \ldots, n \rangle .
\end{gathered} \label{eq7}
\end{equation}
otherwise.
\end{prop}

\dimo If $n$ is odd, an open Heegaard diagram for $M_{25}(n)$ arises from Fig.~\ref{fig:m25n} with the discs corresponding to the faces of $\mathcal P_n$ and the segments of curves being dual to the edges of $\mathcal P_n$. Thus, $\pi_1 (M_{25}(n))$ is generated by $a_1, \ldots, a_n$, $b_1, \ldots, b_n$, $c_1, \ldots, c_n$, $d$. Obviously, passing along the curves of the Heegaard diagram we get defining relations as $({\rm I}_i)$, $({\rm II}_i)$, and $({ \rm III}_i)$ in Proposition~\ref{prop:m25}: $a_i b_{i+1} d^{-1} = 1$,  $a_i c_{i+2}^{-1} b_i^{-1} = 1$,  and $c_i c_{i+1} b_{i-1}^{-1}= 1$, with the additional relation $a_1 a_3 \ldots a_n a_2 \ldots a_{n-1} =1$.

If $n$ is even, the discs corresponding to the faces of $\mathcal P_n$ still give a complete system of meridian discs for a Heegaard surface $F$ of $M_{25}(n)$. The system of curves, defined on $F$, by the edges dual to the edges of $\mathcal P_n$ is proper but not reduced: in fact, by cutting $F$ along these curves, we get two discs. As a consequence, the two systems of curves on $F$ define a generalized Heegaard diagram for $M$ (see \cite{Cattabriga-Mulazzani-Vesnin} for details).

Again this diagram yields a presentation for $\pi_1 (M_{25}(n))$, whose generators are $a_1, \ldots, a_n$, $b_1, \ldots, b_n$, $c_1, \ldots, c_n$, $d$ and whose relations are still $({\rm I}_i)$, $({\rm II}_i)$, and $({ \rm III}_i)$ in Proposition~\ref{prop:m25}: $a_i b_{i+1} d^{-1} = 1$,  $a_i c_{i+2}^{-1} b_i^{-1} = 1$,  and $c_i c_{i+1} b_{i-1}^{-1}= 1$, with two additional relations: $a_1 a_3 \ldots a_{n-1} = 1$ and $a_2 a_4 \ldots a_n =1$. 

Therefore, both for $n$ odd and for $n$ even, the presentations are geometric in the sense that they correspond to (generalized) Heegaard diagrams.
\qed

\begin{cor} \label{cor-m25n-presentation}
The fundamental group of $M_{25}(n)$, $n\geqslant 1$, has presentation:
\begin{equation}
\langle c_1, \ldots, c_n \, | \,   \prod_{j=0}^{k-1} (c_{2+2j} c_{3+2j}^2 )  \prod_{j=0}^{k-1} ( c_{3+2j} c_{4+2j}^2 ) = 1,    \quad
c_{i} c_{i+1}^3 c_{i+2}  = c_{i+1} c_{i+2}^3 c_{i+3}  , \quad i=1, \ldots, n \rangle , \label{eq8}
\end{equation}
if $n=2k+1$, and
\begin{equation}
\langle c_1, \ldots, c_n \, | \,  \prod_{j=0}^{k-1} ( c_{2+2j} c_{3+2j}^2 ) = 1,  \quad  \prod_{j=0}^{k-1} ( c_{3+2j} c_{4+2j}^2 ) = 1,  \quad  c_{i} c_{i+1}^3 c_{i+2}  = c_{i+1} c_{i+2}^3 c_{i+3}   , \quad i=1, \ldots, n \rangle , \label{eq9}
\end{equation}
if $n=2k$, where indices are considered mod $n$.
\end{cor}

\dimo
Let us start with the group presentation given in Proposition~\ref{proposition-m25n-presentation}. Since $b_{i-1} = c_i c_{i+1}$, we get  $a_i = c_{i+1} c_{i+2}^2$ for $i=1, \ldots n$; moreover, $d = c_{i+1} c_{i+2}^3 c_{i+3}$, for all  $i=1,\ldots,n$ (all indices are taken mod~$n$).
\qed

\begin{figure}[h]
\begin{center}
\unitlength=0.7mm
\begin{picture}(180,80)(-20,-5)
\thicklines
\qbezier(-20,0)(-20,0)(100,0) \qbezier(0,60)(0,60)(120,60)
\qbezier(20,0)(20,0)(20,40) \qbezier(60,0)(60,0)(60,40)
\qbezier(100,0)(100,0)(100,40)
\qbezier(0,20)(0,20)(0,60) \qbezier(40,20)(40,20)(40,60)
\qbezier(80,20)(80,20)(80,60) \qbezier(120,20)(120,20)(120,60)
\qbezier(20,40)(20,40)(0,20) \qbezier(20,40)(20,40)(0,60)
\qbezier(20,40)(20,40)(40,20)  \qbezier(20,40)(20,40)(40,60)
\qbezier(60,40)(60,40)(40,20) \qbezier(60,40)(60,40)(40,60)
\qbezier(60,40)(60,40)(80,20)  \qbezier(60,40)(60,40)(80,60)
\qbezier(100,40)(100,40)(80,20) \qbezier(100,40)(100,40)(80,60)
\qbezier(100,40)(100,40)(120,20)  \qbezier(100,40)(100,40)(120,60)
\qbezier(0,20)(0,20)(-20,0)
\qbezier(0,20)(0,20)(20,0) \qbezier(40,20)(40,20)(20,0)
\qbezier(40,20)(40,20)(60,0)  \qbezier(80,20)(80,20)(60,0)
\qbezier(80,20)(80,20)(100,0) \qbezier(120,20)(120,20)(100,0)
\qbezier(120,60)(120,60)(160,60)  \qbezier(160,60)(160,60)(160,20)
\qbezier(100,0)(100,0)(140,0)  \qbezier(140,0)(140,0)(160,20)
\qbezier(120,20)(120,20)(140,0) \qbezier(140,40)(140,40)(140,0)
\qbezier(140,40)(140,40)(120,60) \qbezier(140,40)(140,40)(160,60)
\qbezier(140,40)(140,40)(120,20) \qbezier(140,40)(140,40)(160,20)
\put(140,65){\makebox(0,0)[cc]{$x_1$}}
\put(142,60){\vector(-1,0){4}}
\put(100,65){\makebox(0,0)[cc]{$x_2$}}
\put(102,60){\vector(-1,0){4}}
\put(60,65){\makebox(0,0)[cc]{$x_3$}}
\put(62,60){\vector(-1,0){4}}
\put(20,65){\makebox(0,0)[cc]{$x_4$}}
\put(22,60){\vector(-1,0){4}}
\put(120,-5){\makebox(0,0)[cc]{$x_3$}}
\put(122,0){\vector(-1,0){4}}
\put(80,-5){\makebox(0,0)[cc]{$x_4$}}
\put(82,0){\vector(-1,0){4}}
\put(40,-5){\makebox(0,0)[cc]{$x_1$}}
\put(42,0){\vector(-1,0){4}}
\put(0,-5){\makebox(0,0)[cc]{$x_2$}}
\put(2,0){\vector(-1,0){4}}
\put(165,40){\makebox(0,0)[cc]{$x_3$}}
\put(160,42){\vector(0,-1){4}}
\put(123,40){\makebox(0,0)[cc]{$x_4$}}
\put(120,42){\vector(0,-1){4}}
\put(83,40){\makebox(0,0)[cc]{$x_1$}}
\put(80,42){\vector(0,-1){4}}
\put(43,40){\makebox(0,0)[cc]{$x_2$}}
\put(40,42){\vector(0,-1){4}}
\put(-5,40){\makebox(0,0)[cc]{$x_3$}}
\put(0,42){\vector(0,-1){4}}
\put(132,52){\makebox(0,0)[cc]{$y_1$}}
\put(128,52){\vector(1,-1){4}}
\put(92,52){\makebox(0,0)[cc]{$y_2$}}
\put(88,52){\vector(1,-1){4}}
\put(52,52){\makebox(0,0)[cc]{$y_3$}}
\put(48,52){\vector(1,-1){4}}
\put(12,52){\makebox(0,0)[cc]{$y_4$}}
\put(8,52){\vector(1,-1){4}}
\put(148,52){\makebox(0,0)[cc]{$u$}}
\put(152,52){\vector(-1,-1){4}}
\put(108,52){\makebox(0,0)[cc]{$v$}}
\put(112,52){\vector(-1,-1){4}}
\put(68,52){\makebox(0,0)[cc]{$u$}}
\put(72,52){\vector(-1,-1){4}}
\put(28,52){\makebox(0,0)[cc]{$v$}}
\put(32,52){\vector(-1,-1){4}}
\put(148,12){\makebox(0,0)[cc]{$z_3$}}
\put(152,12){\vector(-1,-1){4}}
\put(108,12){\makebox(0,0)[cc]{$z_4$}}
\put(112,12){\vector(-1,-1){4}}
\put(68,12){\makebox(0,0)[cc]{$z_1$}}
\put(72,12){\vector(-1,-1){4}}
\put(28,12){\makebox(0,0)[cc]{$z_2$}}
\put(32,12){\vector(-1,-1){4}}
\put(-12,12){\makebox(0,0)[cc]{$z_3$}}
\put(-8,12){\vector(-1,-1){4}}
\put(128,32){\makebox(0,0)[cc]{$z_1$}}
\put(132,32){\vector(-1,-1){4}}
\put(88,32){\makebox(0,0)[cc]{$z_2$}}
\put(92,32){\vector(-1,-1){4}}
\put(48,32){\makebox(0,0)[cc]{$z_3$}}
\put(52,32){\vector(-1,-1){4}}
\put(8,32){\makebox(0,0)[cc]{$z_4$}}
\put(12,32){\vector(-1,-1){4}}
\put(152,32){\makebox(0,0)[cc]{$y_3$}}
\put(152,28){\vector(-1,1){4}}
\put(112,32){\makebox(0,0)[cc]{$y_4$}}
\put(112,28){\vector(-1,1){4}}
\put(72,32){\makebox(0,0)[cc]{$y_1$}}
\put(72,28){\vector(-1,1){4}}
\put(32,32){\makebox(0,0)[cc]{$y_2$}}
\put(32,28){\vector(-1,1){4}}
\put(132,12){\makebox(0,0)[cc]{$y_4$}}
\put(132,8){\vector(-1,1){4}}
\put(92,12){\makebox(0,0)[cc]{$y_1$}}
\put(92,8){\vector(-1,1){4}}
\put(52,12){\makebox(0,0)[cc]{$y_2$}}
\put(52,8){\vector(-1,1){4}}
\put(12,12){\makebox(0,0)[cc]{$y_3$}}
\put(12,8){\vector(-1,1){4}}
\put(143,20){\makebox(0,0)[cc]{$z_4$}}
\put(140,18){\vector(0,1){4}}
\put(103,20){\makebox(0,0)[cc]{$z_1$}}
\put(100,18){\vector(0,1){4}}
\put(63,20){\makebox(0,0)[cc]{$z_2$}}
\put(60,18){\vector(0,1){4}}
\put(23,20){\makebox(0,0)[cc]{$z_3$}}
\put(20,18){\vector(0,1){4}}
\put(140,55){\makebox(0,0)[cc]{$\bf A_1$}}
\put(100,55){\makebox(0,0)[cc]{$\bf A_2$}}
\put(60,55){\makebox(0,0)[cc]{$\bf A_3$}}
\put(20,55){\makebox(0,0)[cc]{$\bf A_4$}}
\put(150,40){\makebox(0,0)[cc]{$\bf \bar{A}_3$}}
\put(130,40){\makebox(0,0)[cc]{$\bf B_1$}}
\put(110,40){\makebox(0,0)[cc]{$\bf \bar{A}_4$}}
\put(90,40){\makebox(0,0)[cc]{$\bf B_2$}}
\put(70,40){\makebox(0,0)[cc]{$\bf \bar{A}_1$}}
\put(50,40){\makebox(0,0)[cc]{$\bf B_3$}}
\put(30,40){\makebox(0,0)[cc]{$\bf \bar{A}_2$}}
\put(10,40){\makebox(0,0)[cc]{$\bf B_4$}}
\put(150,20){\makebox(0,0)[cc]{$\bf \bar{C}_1$}}
\put(130,20){\makebox(0,0)[cc]{$\bf C_2$}}
\put(110,20){\makebox(0,0)[cc]{$\bf \bar{C}_2$}}
\put(90,20){\makebox(0,0)[cc]{$\bf C_3$}}
\put(70,20){\makebox(0,0)[cc]{$\bf \bar{C}_3$}}
\put(50,20){\makebox(0,0)[cc]{$\bf C_4$}}
\put(30,20){\makebox(0,0)[cc]{$\bf \bar{C}_4$}}
\put(10,20){\makebox(0,0)[cc]{$\bf C_1$}}
\put(120,5){\makebox(0,0)[cc]{$\bf \bar{B}_4$}}
\put(80,5){\makebox(0,0)[cc]{$\bf \bar{B}_1$}}
\put(40,5){\makebox(0,0)[cc]{$\bf \bar{B}_2$}}
\put(0,5){\makebox(0,0)[cc]{$\bf \bar{B}_3$}}
\put(70,70){\makebox(0,0)[cc]{$\bf D$}}
\put(70,-10){\makebox(0,0)[cc]{$\bf \bar{D}$}}
\put(160,60){\makebox(0,0)[cc]{$\blacklozenge$}}
\put(120,60){\makebox(0,0)[cc]{$\blacksquare$}}
\put(80,60){\makebox(0,0)[cc]{$\blacklozenge$}}
\put(40,60){\makebox(0,0)[cc]{$\blacksquare$}}
\put(0,60){\makebox(0,0)[cc]{$\blacklozenge$}}
\put(140,40){\makebox(0,0)[cc]{$\blacksquare$}}
\put(100,40){\makebox(0,0)[cc]{$\blacklozenge$}}
\put(60,40){\makebox(0,0)[cc]{$\blacksquare$}}
\put(20,40){\makebox(0,0)[cc]{$\blacklozenge$}}
\put(160,20){\makebox(0,0)[cc]{$\blacksquare$}}
\put(120,20){\makebox(0,0)[cc]{$\blacklozenge$}}
\put(80,20){\makebox(0,0)[cc]{$\blacksquare$}}
\put(40,20){\makebox(0,0)[cc]{$\blacklozenge$}}
\put(0,20){\makebox(0,0)[cc]{$\blacksquare$}}
\put(140,0){\makebox(0,0)[cc]{$\blacklozenge$}}
\put(100,0){\makebox(0,0)[cc]{$\blacksquare$}}
\put(60,0){\makebox(0,0)[cc]{$\blacklozenge$}}
\put(20,0){\makebox(0,0)[cc]{$\blacksquare$}}
\put(-20,0){\makebox(0,0)[cc]{$\blacklozenge$}}
\end{picture}
\end{center} \caption{Equivalence of faces in the construction of $M_{25}(4)$.} \label{fig:m25-4}
\end{figure}
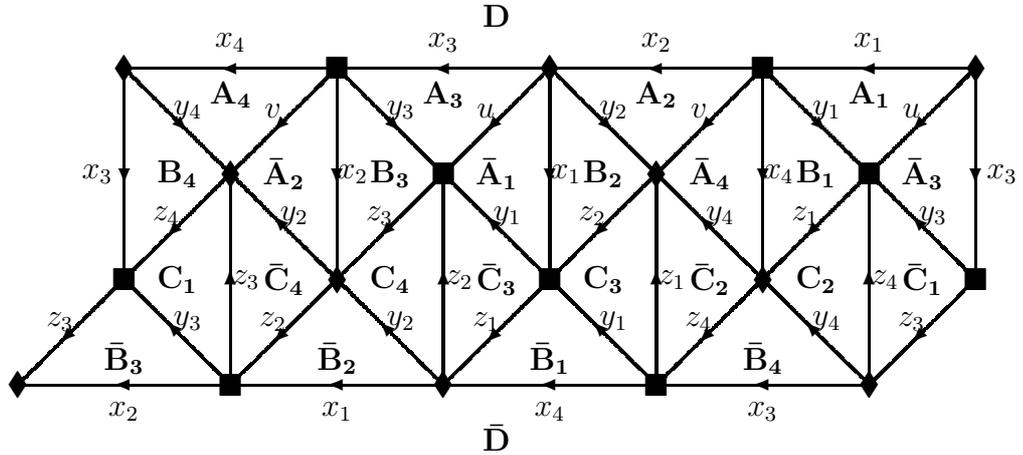

It is stated in~\cite[p.~391]{Cavicchioli2} that for $n \geqslant 3$ the group $G(n)$ with the presentation (\ref{group-G}),  where the relations correspond to the boundaries of the 2-faces, as shown for $M_{25}(3)$ in  Fig.~\ref{fig:m25}, is the fundamental group of $M_{25}(n)$. We observe that $G(n)$ is isomorphic to $\pi_1 (M_{25}(n))$ only if $n$ is odd. This is the case when the quotient space $\mathcal P_n / \psi_n$ has exactly one vertex. If there is more than one vertex, more careful considerations are necessary (see, for example, \cite[Section~62]{Seifert-Threlfall}). If $n = 2k$, then, as was already pointed out above, all vertices of $\mathcal P_{n}$ will form two classes of equivalence (see Fig.~\ref{fig:m25-4} for $M_{25}(4)$ where one class of vertices is marked by $\blacksquare$ and another by $\blacklozenge$) and all its edges will form $(3n+2)$ classes of equivalence.
\noindent
For example, as we can see from the polyhedral schemata for $M_{25}(4)$ in Fig.~\ref{fig:m25-4}, there are two edges in the class $u$, and two edges in the class $v$; on the contrary, from the presentation (\ref{group-G}), four edges in the class $u$ would be expected.  Actually, for $n$ even, we must add a relation corresponding to a maximal tree of the 1-skeleton of the complex, i.e. an edge. For instance we suppose $v=1$: indeed $v$ is an edge connecting a vertex from class $\blacksquare$ with a vertex from class $\blacklozenge$.
Therefore, for $n = 2k$,  the fundamental group $\pi_1 (M_{25}(n))$ is isomorphic to $H(n)$ with the following presentation:
\begin{equation}
\begin{gathered}
H(n) = \langle x_1, \ldots, x_n;  y_1, \ldots, y_n; z_1, \ldots, z_n;  u \quad | \quad x_1 x_2 \ldots x_n = 1, \\
y_i z_i  = x_{i-1}, \quad z_{i-1} z_i  = y_{i-1} , \quad i=1, \ldots, n , \qquad  x_{2j-1} y_{2j-1} = u, \quad x_{2j} y_{2j} = 1 , \quad j = 1, \ldots k \rangle .
\end{gathered} \label{correct25}
\end{equation}

\section{Covering properties and volumes of manifolds $M_{25}(n)$}
\label{sec4}

\begin{prop} \label{prop-M25-2}
The manifold $M_{25}(2)$ is the lens space $L_{3,1}$.
\end{prop}

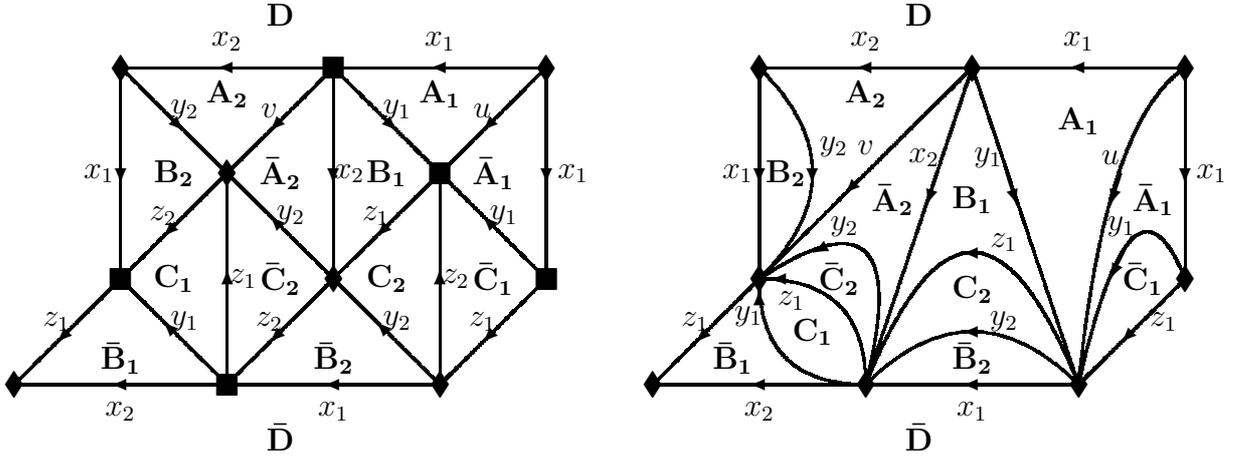
\begin{figure}[h]
\begin{center}
\unitlength=0.7mm
\begin{picture}(180,80)
\put(0,0){\begin{picture}(180,80)(80,-5)
\thicklines
\qbezier(60,0)(60,0)(100,0) \qbezier(80,60)(80,60)(120,60)
\qbezier(80,20)(80,20)(80,60) \qbezier(120,20)(120,20)(120,60)
\qbezier(100,40)(100,40)(80,20) \qbezier(100,40)(100,40)(80,60)
\qbezier(100,40)(100,40)(120,20)  \qbezier(100,40)(100,40)(120,60)
\qbezier(80,20)(80,20)(60,0) \qbezier(100,40)(100,40)(100,0)
\qbezier(80,20)(80,20)(100,0) \qbezier(120,20)(120,20)(100,0)
\qbezier(120,60)(120,60)(160,60)  \qbezier(160,60)(160,60)(160,20)
\qbezier(100,0)(100,0)(140,0)  \qbezier(140,0)(140,0)(160,20)
\qbezier(120,20)(120,20)(140,0) \qbezier(140,40)(140,40)(140,0)
\qbezier(140,40)(140,40)(120,60) \qbezier(140,40)(140,40)(160,60)
\qbezier(140,40)(140,40)(120,20) \qbezier(140,40)(140,40)(160,20)
\put(140,65){\makebox(0,0)[cc]{$x_1$}}
\put(142,60){\vector(-1,0){4}}
\put(100,65){\makebox(0,0)[cc]{$x_2$}}
\put(102,60){\vector(-1,0){4}}
\put(120,-5){\makebox(0,0)[cc]{$x_1$}}
\put(122,0){\vector(-1,0){4}}
\put(80,-5){\makebox(0,0)[cc]{$x_2$}}
\put(82,0){\vector(-1,0){4}}
\put(165,40){\makebox(0,0)[cc]{$x_1$}}
\put(160,42){\vector(0,-1){4}}
\put(123,40){\makebox(0,0)[cc]{$x_2$}}
\put(120,42){\vector(0,-1){4}}
\put(76,40){\makebox(0,0)[cc]{$x_1$}}
\put(80,42){\vector(0,-1){4}}
\put(132,52){\makebox(0,0)[cc]{$y_1$}}
\put(128,52){\vector(1,-1){4}}
\put(92,52){\makebox(0,0)[cc]{$y_2$}}
\put(88,52){\vector(1,-1){4}}
\put(148,52){\makebox(0,0)[cc]{$u$}}
\put(152,52){\vector(-1,-1){4}}
\put(108,52){\makebox(0,0)[cc]{$v$}}
\put(112,52){\vector(-1,-1){4}}
\put(148,12){\makebox(0,0)[cc]{$z_1$}}
\put(152,12){\vector(-1,-1){4}}
\put(108,12){\makebox(0,0)[cc]{$z_2$}}
\put(112,12){\vector(-1,-1){4}}
\put(68,12){\makebox(0,0)[cc]{$z_1$}}
\put(72,12){\vector(-1,-1){4}}
\put(128,32){\makebox(0,0)[cc]{$z_1$}}
\put(132,32){\vector(-1,-1){4}}
\put(88,32){\makebox(0,0)[cc]{$z_2$}}
\put(92,32){\vector(-1,-1){4}}
\put(152,32){\makebox(0,0)[cc]{$y_1$}}
\put(152,28){\vector(-1,1){4}}
\put(112,32){\makebox(0,0)[cc]{$y_2$}}
\put(112,28){\vector(-1,1){4}}
\put(132,12){\makebox(0,0)[cc]{$y_2$}}
\put(132,8){\vector(-1,1){4}}
\put(92,12){\makebox(0,0)[cc]{$y_1$}}
\put(92,8){\vector(-1,1){4}}
\put(143,20){\makebox(0,0)[cc]{$z_2$}}
\put(140,18){\vector(0,1){4}}
\put(103,20){\makebox(0,0)[cc]{$z_1$}}
\put(100,18){\vector(0,1){4}}
\put(140,55){\makebox(0,0)[cc]{$\bf A_1$}}
\put(100,55){\makebox(0,0)[cc]{$\bf A_2$}}
\put(150,40){\makebox(0,0)[cc]{$\bf \bar{A}_1$}}
\put(130,40){\makebox(0,0)[cc]{$\bf B_1$}}
\put(110,40){\makebox(0,0)[cc]{$\bf \bar{A}_2$}}
\put(90,40){\makebox(0,0)[cc]{$\bf B_2$}}
\put(150,20){\makebox(0,0)[cc]{$\bf \bar{C}_1$}}
\put(130,20){\makebox(0,0)[cc]{$\bf C_2$}}
\put(110,20){\makebox(0,0)[cc]{$\bf \bar{C}_2$}}
\put(90,20){\makebox(0,0)[cc]{$\bf C_1$}}
\put(120,5){\makebox(0,0)[cc]{$\bf \bar{B}_2$}}
\put(80,5){\makebox(0,0)[cc]{$\bf \bar{B}_1$}}
\put(110,70){\makebox(0,0)[cc]{$\bf D$}}
\put(110,-10){\makebox(0,0)[cc]{$\bf \bar{D}$}}
\put(160,60){\makebox(0,0)[cc]{$\blacklozenge$}}
\put(120,60){\makebox(0,0)[cc]{$\blacksquare$}}
\put(80,60){\makebox(0,0)[cc]{$\blacklozenge$}}
\put(140,40){\makebox(0,0)[cc]{$\blacksquare$}}
\put(100,40){\makebox(0,0)[cc]{$\blacklozenge$}}
\put(160,20){\makebox(0,0)[cc]{$\blacksquare$}}
\put(120,20){\makebox(0,0)[cc]{$\blacklozenge$}}
\put(80,20){\makebox(0,0)[cc]{$\blacksquare$}}
\put(140,0){\makebox(0,0)[cc]{$\blacklozenge$}}
\put(100,0){\makebox(0,0)[cc]{$\blacksquare$}}
\put(60,0){\makebox(0,0)[cc]{$\blacklozenge$}}
\end{picture}
}
\put(120,0){\begin{picture}(180,80)(80,-5)
\thicklines
\qbezier(60,0)(60,0)(140,0) \qbezier(80,60)(80,60)(160,60)
\qbezier(80,20)(80,20)(80,60)  \qbezier(80,20)(80,20)(60,0)
 \qbezier(160,60)(160,60)(160,20) \qbezier(140,0)(140,0)(160,20)
\qbezier(80,60)(100,40)(80,20)
\put(94,45){\makebox(0,0)[cc]{$y_{2}$}}
\put(90,42){\vector(0,-1){4}}
\put(85,40){\makebox(0,0)[cc]{$\bf B_2$}}
\qbezier(120,60)(120,60)(80,20)
\put(100,45){\makebox(0,0)[cc]{$v$}}
\put(100,40){\vector(-1,-1){4}}
\put(100,55){\makebox(0,0)[cc]{$\bf A_2$}}
\qbezier(120,60)(120,60)(100,0)
\put(111,43){\makebox(0,0)[cc]{$x_{2}$}}
\put(113.5,40){\vector(-1,-3){2}}
\put(105,35){\makebox(0,0)[cc]{$\bf \bar{A}_2$}}
\qbezier(120,60)(120,60)(140,0)
\put(123,43){\makebox(0,0)[cc]{$y_{1}$}}
\put(126.5,40){\vector(1,-3){2}}
\put(120,35){\makebox(0,0)[cc]{$\bf B_1$}}
\qbezier(160,60)(145,50)(140,0)
\put(146,43){\makebox(0,0)[cc]{$u$}}
\put(148,40){\vector(-1,-3){2}}
\put(154,35){\makebox(0,0)[cc]{$\bf \bar{A}_1$}}
\qbezier(80,20)(80,0)(100,0)
\put(75,5){\makebox(0,0)[cc]{$\bf \bar{B}_1$}}
\put(78,13){\makebox(0,0)[cc]{$y_1$}}
\put(80,16){\vector(0,1){2}}
\qbezier(80,20)(100,20)(100,0)
\put(90,10){\makebox(0,0)[cc]{$\bf C_1$}}
\put(86,16){\makebox(0,0)[cc]{$z_1$}}
\put(84,20){\vector(-1,0){2}}
\qbezier(80,20)(110,40)(100,0)
\put(95,20){\makebox(0,0)[cc]{$\bf \bar{C}_2$}}
\put(96,30){\makebox(0,0)[cc]{$y_2$}}
\put(92,26){\vector(-2,-1){2}}
\qbezier(100,0)(120,20)(140,0)
\put(126,12){\makebox(0,0)[cc]{$y_2$}}
\put(122,10){\vector(-1,0){4}}
\put(120,5){\makebox(0,0)[cc]{$\bf \bar{B}_2$}}
\qbezier(100,0)(120,50)(140,0)
\put(126,27){\makebox(0,0)[cc]{$z_1$}}
\put(122,25){\vector(-1,0){4}}
\put(120,18){\makebox(0,0)[cc]{$\bf C_2$}}
\qbezier(140,0)(150,45)(160,20)
\put(148,30){\makebox(0,0)[cc]{$y_1$}}
\put(150,28){\vector(-1,-2){4}}
\put(152,20){\makebox(0,0)[cc]{$\bf \bar{C}_1$}}
\put(140,65){\makebox(0,0)[cc]{$x_1$}}
\put(142,60){\vector(-1,0){4}}
\put(100,65){\makebox(0,0)[cc]{$x_2$}}
\put(102,60){\vector(-1,0){4}}
\put(120,-5){\makebox(0,0)[cc]{$x_1$}}
\put(122,0){\vector(-1,0){4}}
\put(80,-5){\makebox(0,0)[cc]{$x_2$}}
\put(82,0){\vector(-1,0){4}}
\put(165,40){\makebox(0,0)[cc]{$x_1$}}
\put(160,42){\vector(0,-1){4}}
\put(76,40){\makebox(0,0)[cc]{$x_1$}}
\put(80,42){\vector(0,-1){4}}
\put(156,12){\makebox(0,0)[cc]{$z_1$}}
\put(152,12){\vector(-1,-1){4}}
\put(68,12){\makebox(0,0)[cc]{$z_1$}}
\put(72,12){\vector(-1,-1){4}}
\put(140,50){\makebox(0,0)[cc]{$\bf A_1$}}
\put(110,70){\makebox(0,0)[cc]{$\bf D$}}
\put(110,-10){\makebox(0,0)[cc]{$\bf \bar{D}$}}
\put(160,60){\makebox(0,0)[cc]{$\blacklozenge$}}
\put(120,60){\makebox(0,0)[cc]{$\blacklozenge$}}
\put(80,60){\makebox(0,0)[cc]{$\blacklozenge$}}
\put(160,20){\makebox(0,0)[cc]{$\blacklozenge$}}
\put(80,20){\makebox(0,0)[cc]{$\blacklozenge$}}
\put(140,0){\makebox(0,0)[cc]{$\blacklozenge$}}
\put(100,0){\makebox(0,0)[cc]{$\blacklozenge$}}
\put(60,0){\makebox(0,0)[cc]{$\blacklozenge$}}
\end{picture}
}
\end{picture}
\end{center} \caption{Two-vertex and one-vertex fundamental polyhedra for $M_{25}(2)$.} \label{fig:m25-2}
\end{figure}

\begin{figure}[h]
\centering{
\includegraphics[height=6.cm]{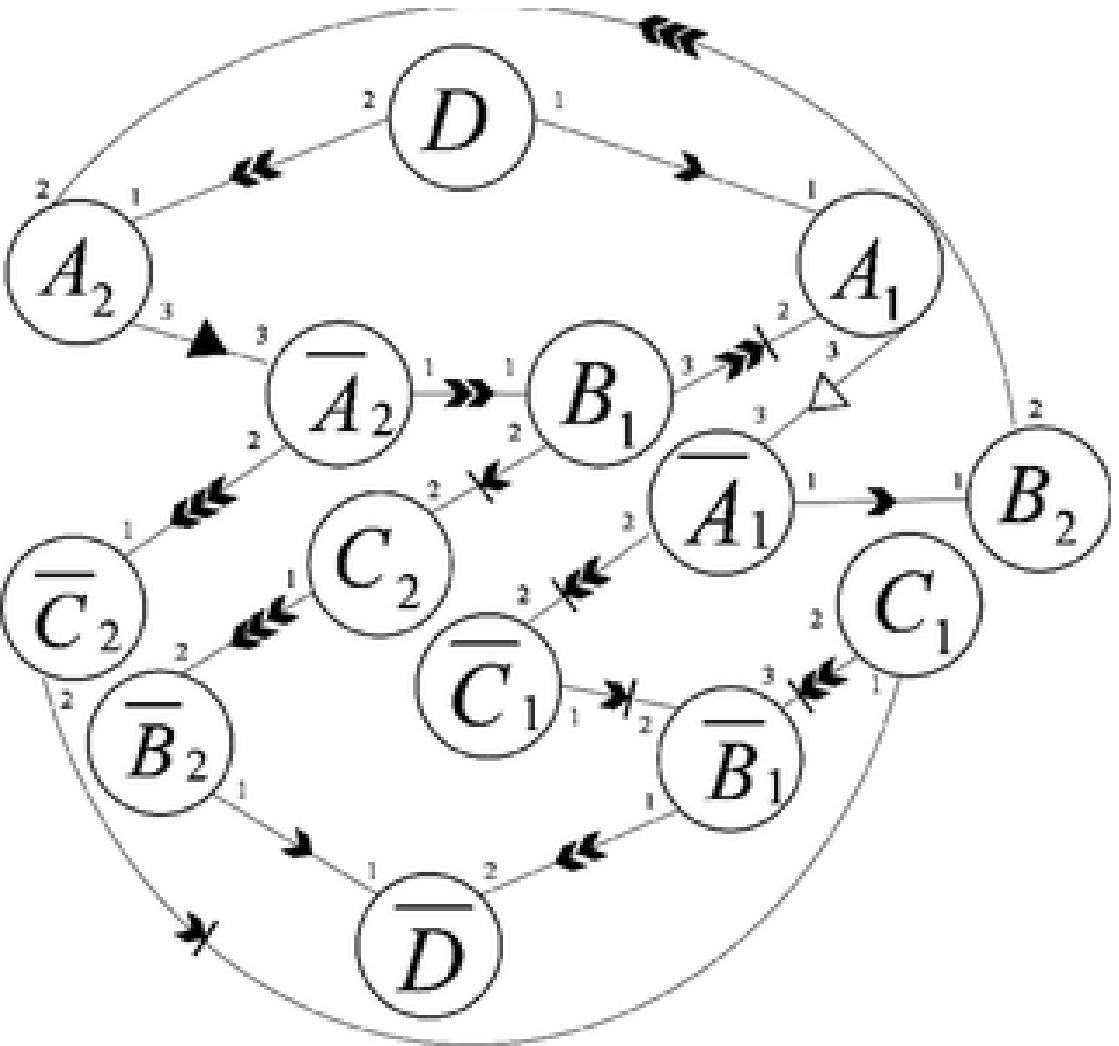}  \qquad \qquad
\includegraphics[height=6.cm]{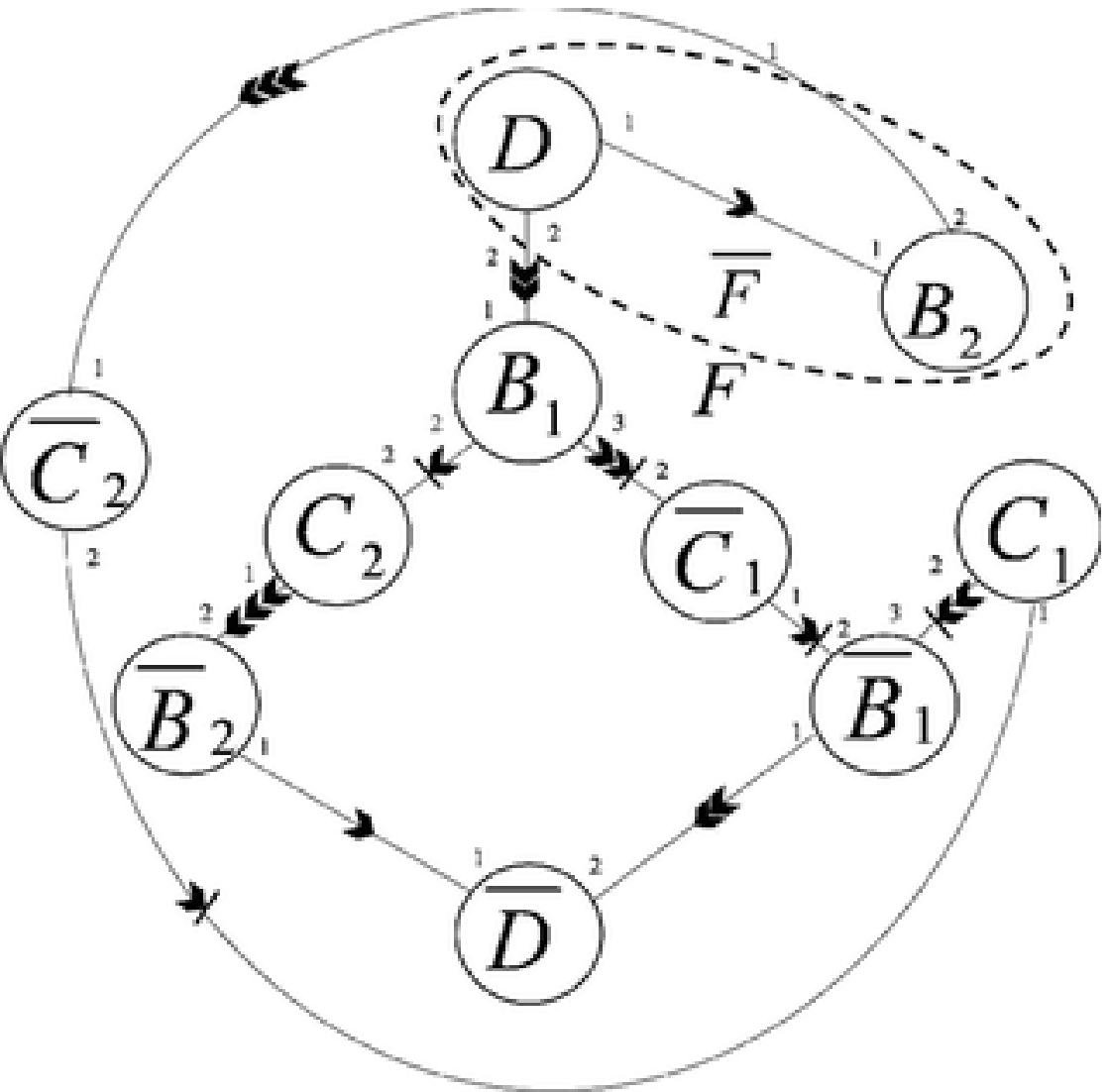}  \\
\includegraphics[height=4.5cm]{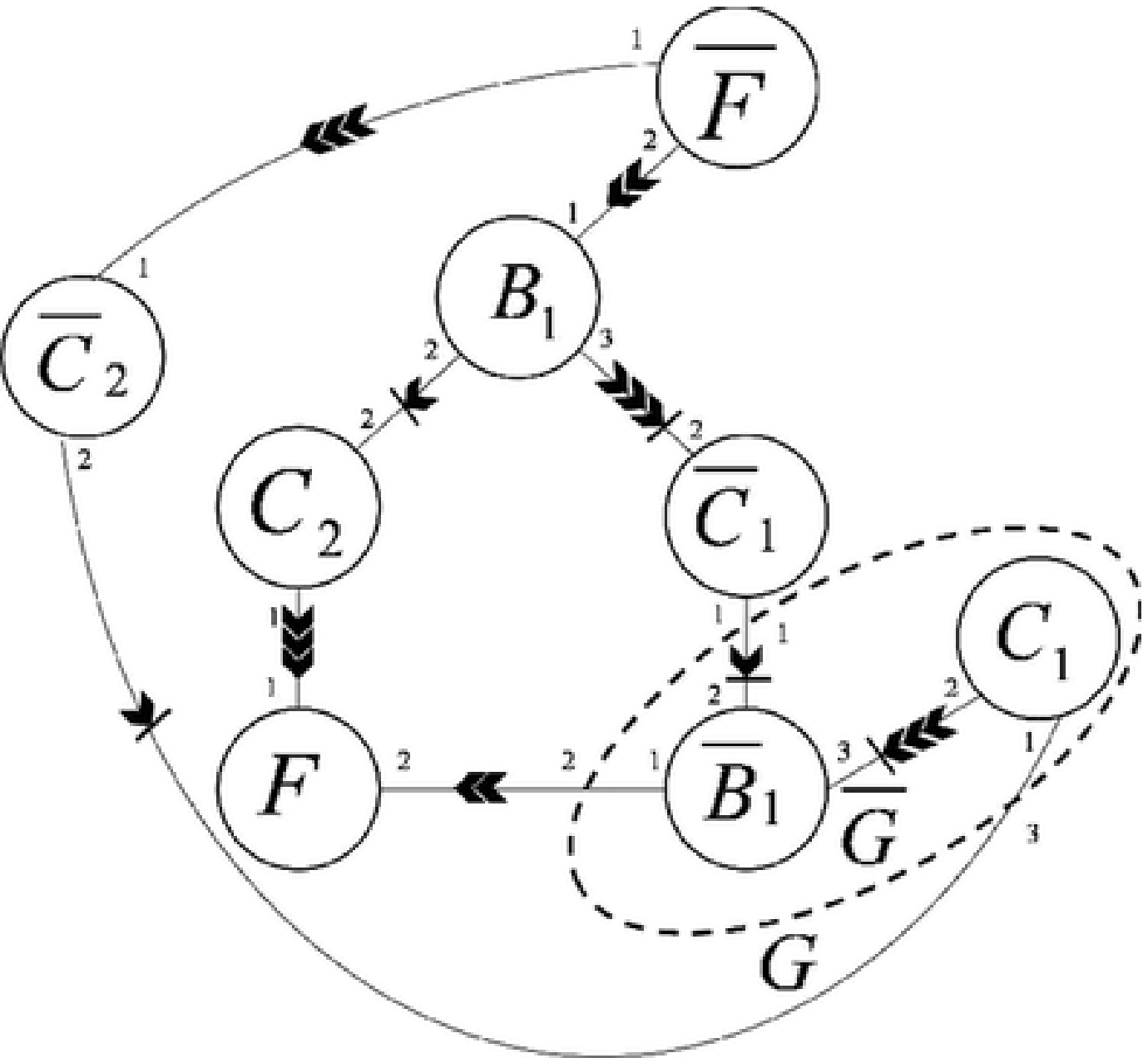}  \qquad \qquad
\includegraphics[height=4.5cm]{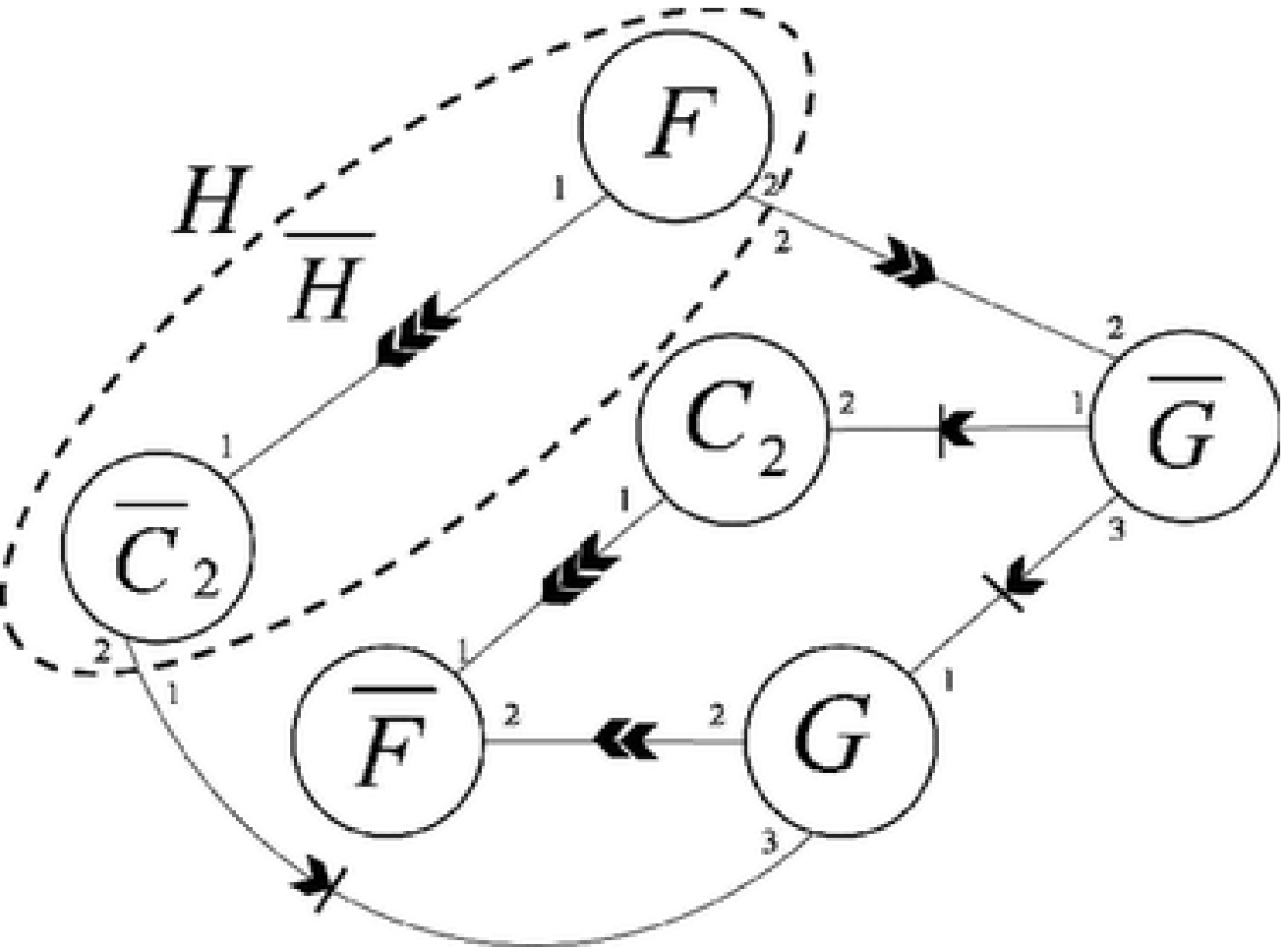}  \\
\includegraphics[height=3.5cm]{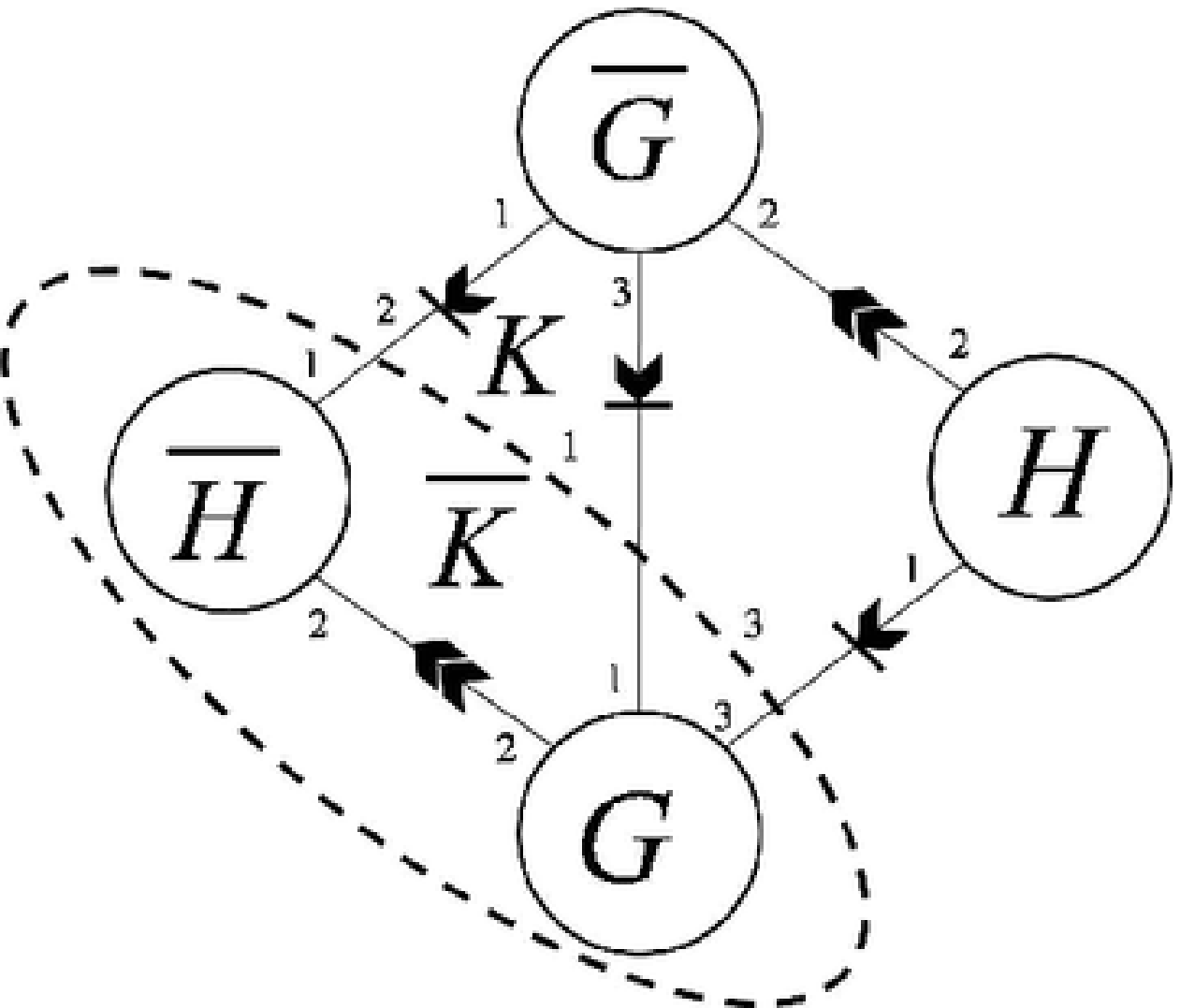} \qquad
\includegraphics[height=3.5cm]{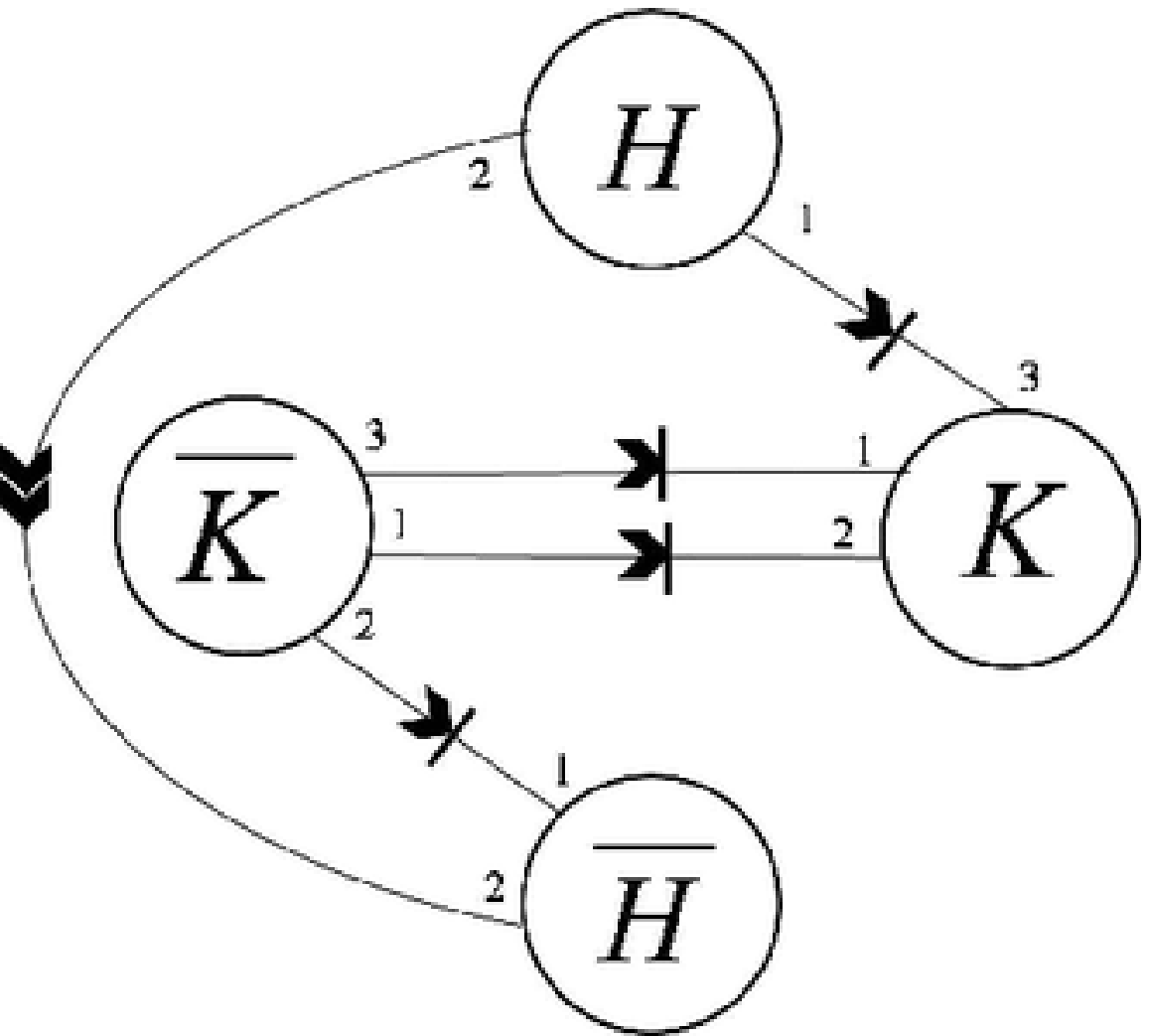}  \qquad
\includegraphics[height=1.cm]{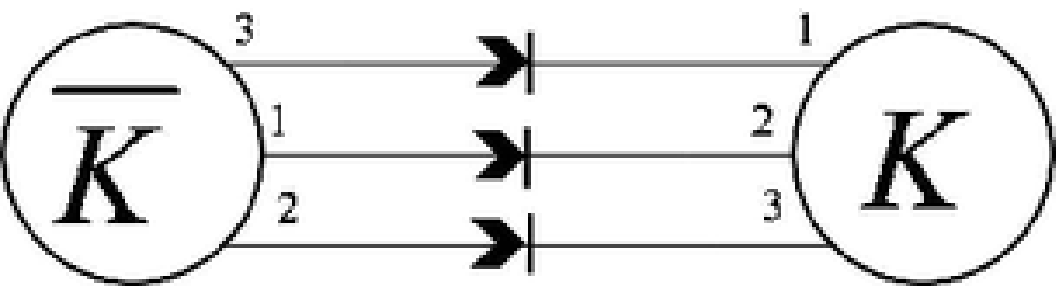}
}
\caption{Heegaard diagrams from $M_{25}(2)$ to $L_{3,1}$.} 
\label{fig:lens_from_m25(2)}
\end{figure}

\dimo
Let us contract the edges of  $\mathcal P_{2}$ labelled by $z_{2}$ to deform it to a fundamental polyhedron whose quotient by $\psi_{2}$ has one vertex (see Fig.~\ref{fig:m25-2}).
We are using the dual on this new polyhedra to get a Heegaard diagram for $M_{25}(2)$. The transformations of Heegaard diagrams from $M_{25}(2)$ to $L_{3,1}$ are drawn in Fig.~\ref{fig:lens_from_m25(2)}.
\qed

With regard to the covering properties of $M_{25}(n)$, we must distinguish again the case $n$ even from that of $n$ odd.

\begin{thm} \label{theorem:m25n-covering}
(1) For every $n$ odd, $n \geqslant 3$, the manifold $M_{25}(n)$ is an $n$-fold strongly-cyclic  branched covering of the lens space $L_{3,1}$, branched over a 2-component link. Moreover, $M_{25}(1)$ is the lens space $L_{3,1}$.

(2) For every $n = 2k$, $k \geqslant 2$, the manifold $M_{25}(n)$ is a $n/2$-fold strongly-cyclic  branched covering of the lens space $L_{3,1}$, branched over a 3-component link.
\end{thm}

\dimo
(1) Suppose that $n$ is odd. Denote by $\rho_n$ the rotational symmetry of $\mathcal P_n$ sending $X_i$ to $X_{i+1}$ with indices taken mod $n$, $i=1, \ldots, n$, where $X$ belongs to the set of letters used for the notations of vertices: $\{ P, Q, R, S\}$. This symmetry induces a cyclic symmetry of the quotient space $M_{25}(n) = \mathcal P_n / \psi_n$, and we denote it by $\rho_n$, too. The quotient space $M_{25}(n) / \rho_n$ is an orbifold whose underlying manifold is $M_{25}(1)$ with a 2-component singular set $\mathcal L_{odd} = \ell_1 \cup \ell_2$. According to the description of the equivalence classes of edges given in Prop.~\ref{prop:m25}, the first component $\ell_1$ of  $\mathcal L_{odd}$ corresponds to the class of edges $\{ P_1 Q_1, \ldots, P_n Q_n\}$. Only one element from this class will appear in the quotient space, so $\ell_1$ has singularity index $n$. The second singular curve $\ell_2$, also having singularity index $n$, corresponds to the axis of rotation $\rho$. Through the Heegaard diagram of $M_{25}(n) / \rho_n$ we understand the Heegaard diagram of $M_{25}(1)$ with information about the singular set $\mathcal L_{odd}$ presented.

\begin{figure}[h]
\centering{
\includegraphics[height=6cm]{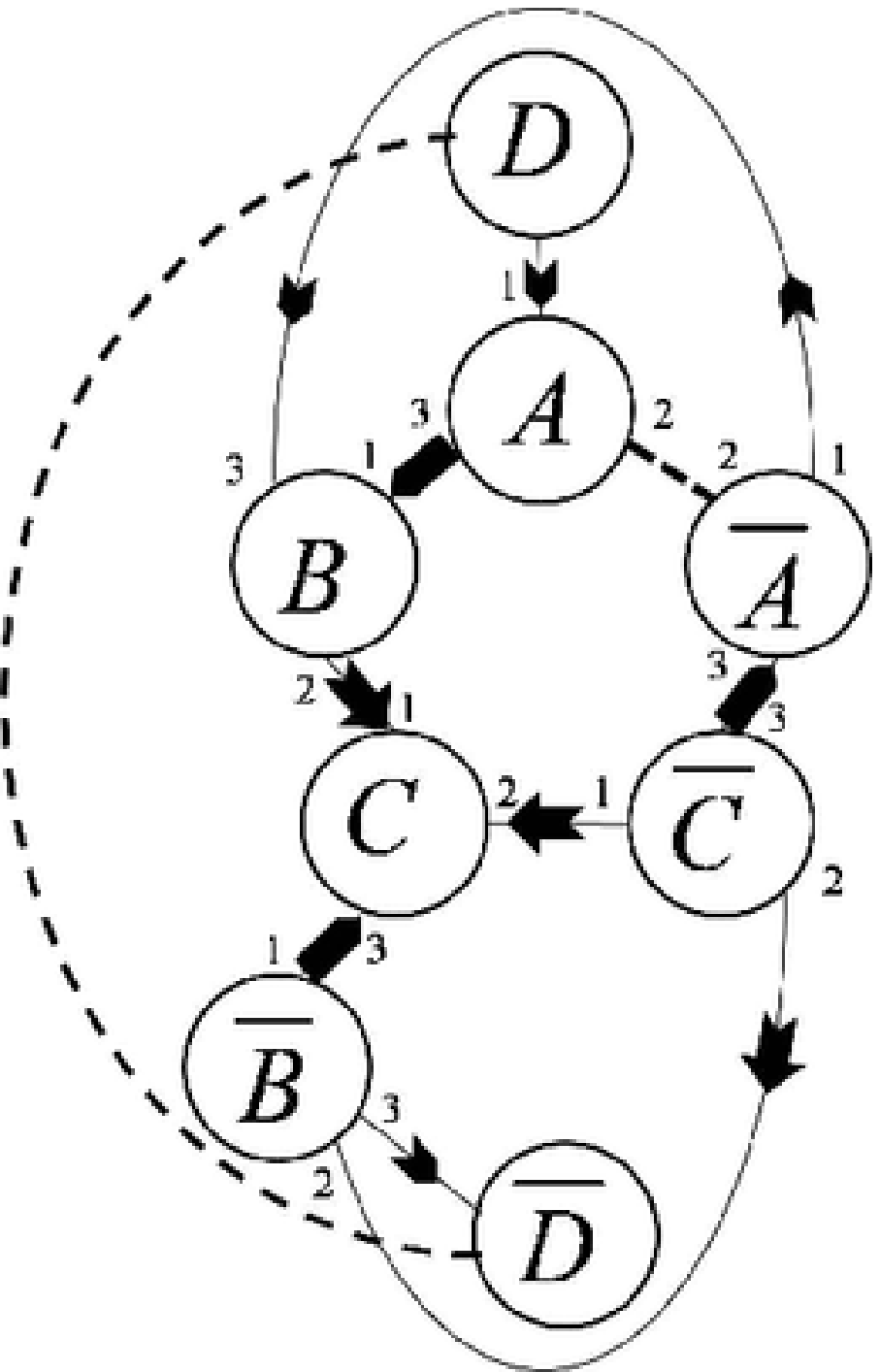}  \qquad
\includegraphics[height=6cm]{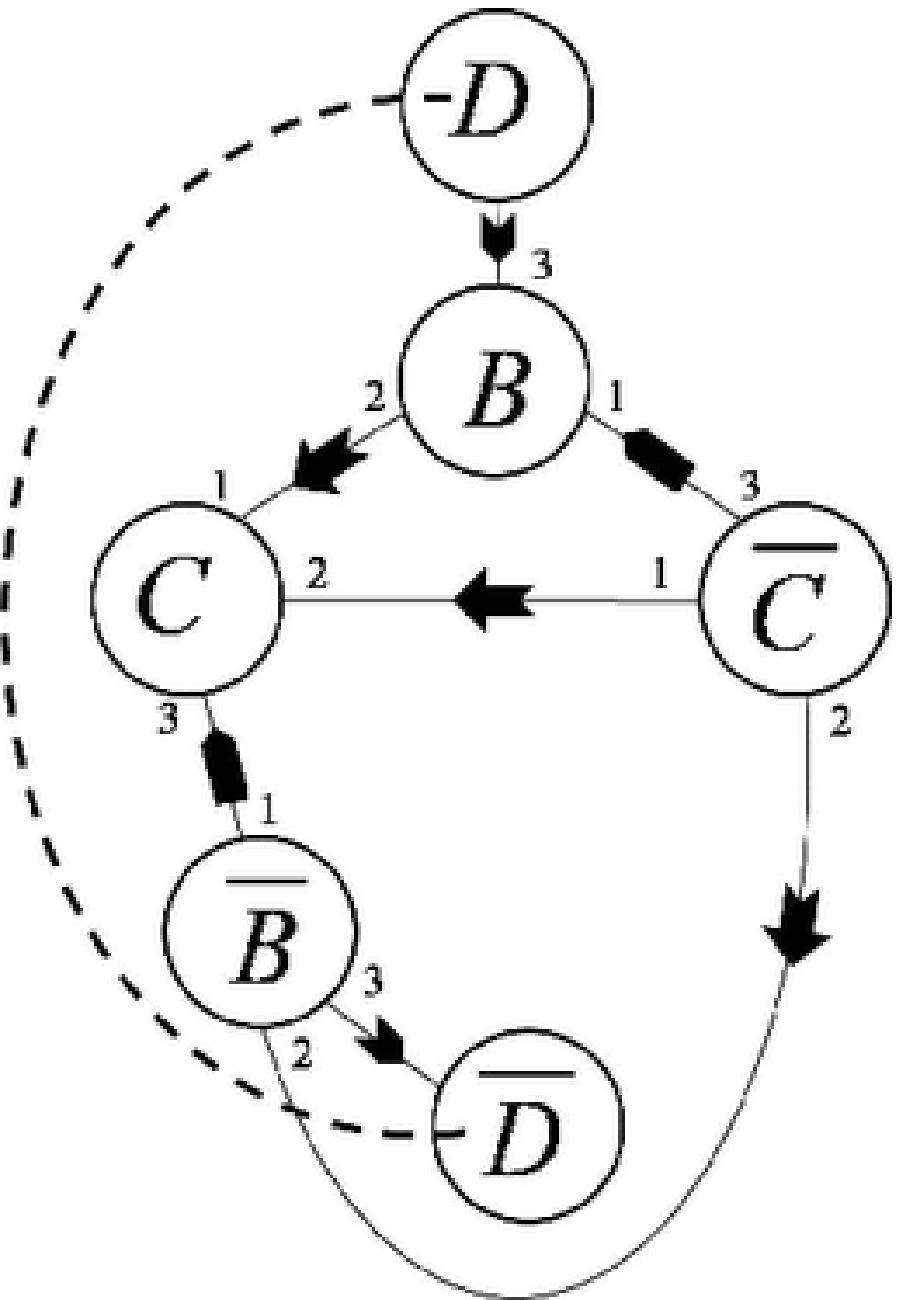} \qquad
\includegraphics[height=4cm]{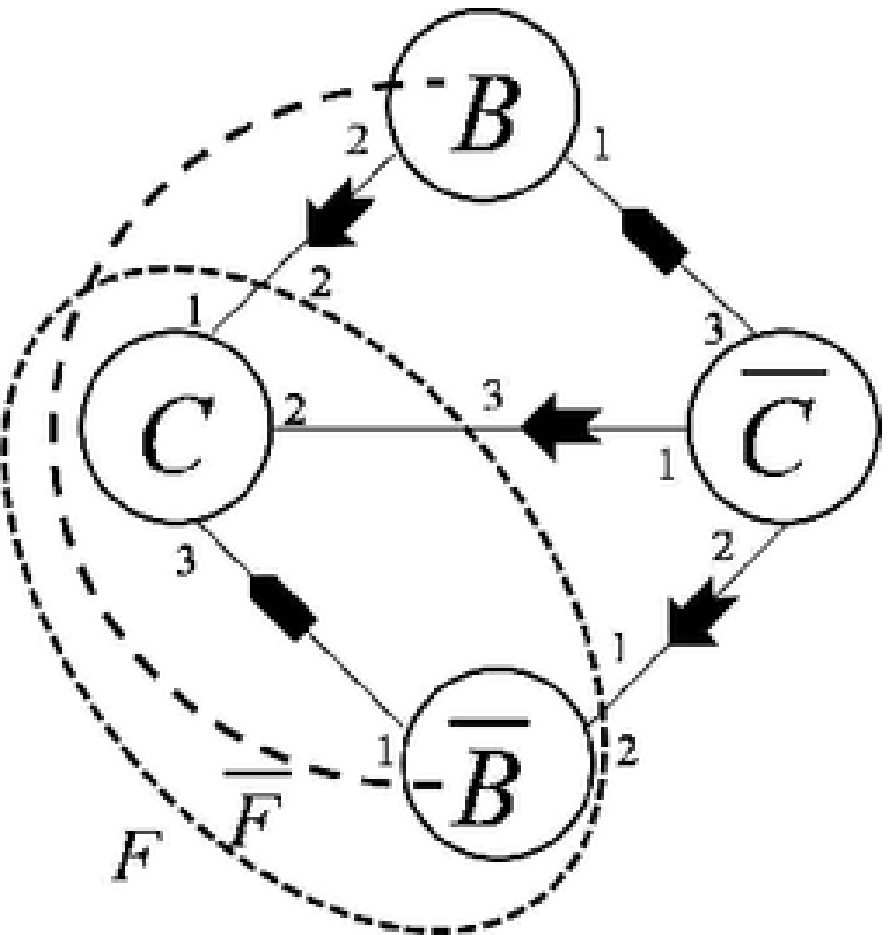}  \\
\includegraphics[height=4cm]{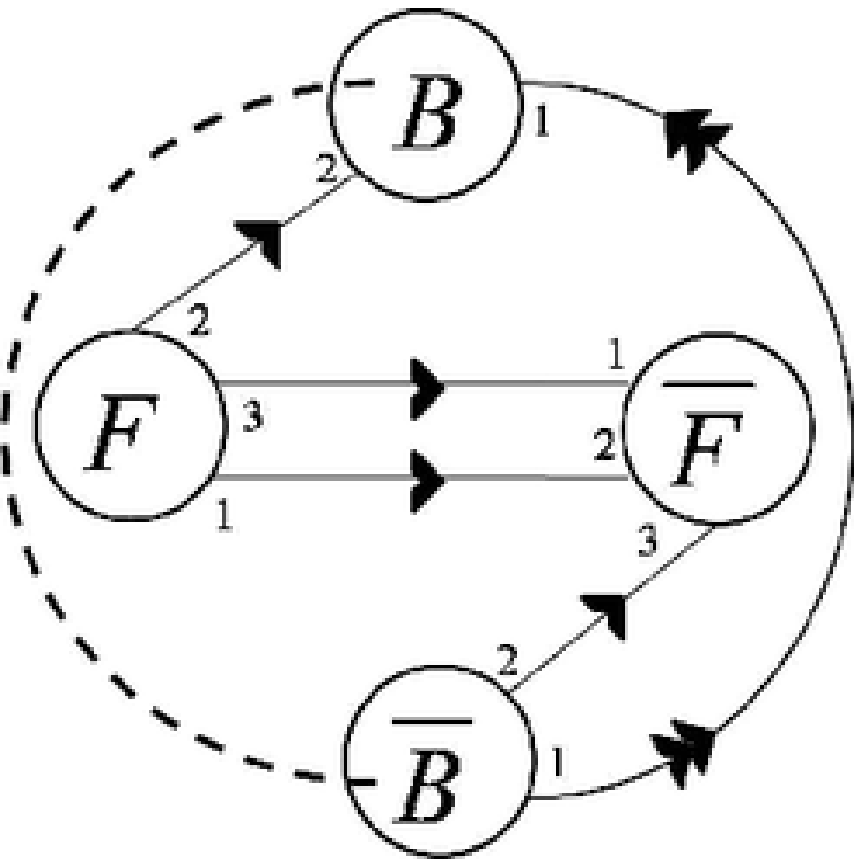}  \qquad \qquad
\includegraphics[height=1.5cm]{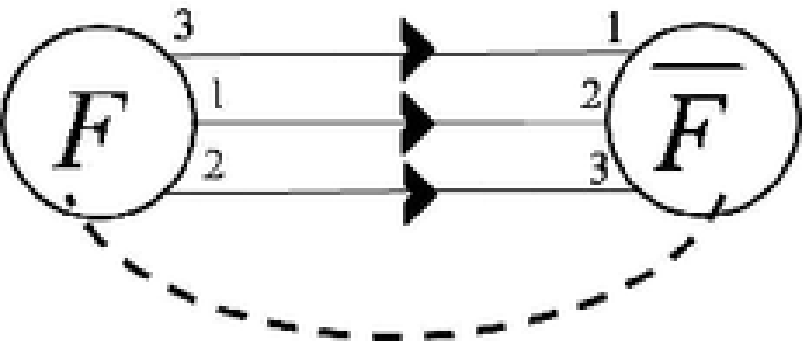}
}
\caption{Heegaard diagrams from $M_{25}(n) / \rho_n$ to $L_{3,1}$.} \label{fig:lens_from_m25}
\end{figure}

The transformations of Heegaard diagrams from $M_{25}(n) / \rho_n$ to $L_{3,1}$ are drawn in Fig.~\ref{fig:lens_from_m25}. Here the first component of the singular set is represented by a dashed segment connecting the discs $A$ and $\bar{A}$; the second component by a dashed segment, connecting the  discs $D$ and $\bar{D}$. Both components have branching index $n$, thus the cyclic covering is strongly cyclic.
At the first step we identify discs  $A$ and $\bar{A}$, forming an 1-handle. The dashed segment (which corresponds to the first component of the singular set) connecting these discs will give a 2-handle to glue up this 1-handle (i.e. they form a pair of complementary handles).  Thus, the first component is a trivial knot; we are not drawing it in next figures.  At the second step we  cancel the discs $D$ and $\bar{D}$, since they are connected only with the discs $B$ and $\bar{B}$, respectively, and the connecting segments are glued together. At the third step we cut along the curve represented by the dotted line to form a new pair of discs $F$ and $\bar{F}$ and then identify the discs $C$ and $\bar{C}$. After that, we easily get a genus one Heegaard diagram, with the discs $F$ and $\bar{F}$, which is the standard diagram for the lens space $L_{3,1}$, where the dotted line represents the second component of the singular set $\mathcal L_{odd}$.

(2) Suppose $n = 2k$, $k \geqslant 2$. Denote by $\rho_k$ the rotational symmetry of $\mathcal P_n$ sending $X_i$ to $X_{i+2}$, $i=1, \ldots, n$, with indices taken mod $n$, where $X$ belongs to the set of letters used for the notations of vertices: $\{ P, Q, R, S\}$. This symmetry induces a cyclic symmetry of the quotient space $M_{25}(n) = \mathcal P_n / \psi_n$, and we denote it by $\rho_k$, too. The quotient space $M_{25}(n) / \rho_k$ is an orbifold whose underlying manifold is $M_{25}(2)$. We observe that the singular set $\mathcal L_{even}$ of the orbifold $M_{25}(n) / \rho_k$ has three components: $\mathcal L_{even} = \ell'_1 \cup \ell'_2 \cup \ell'_3$. Indeed, according to the description of the equivalence classes of the edges given in Prop.~\ref{prop:m25} for the case of $n$ even, the first component $\ell'_1$ corresponds to the class of edges $\{ P_1 Q_1, P_3 Q_3, \ldots, P_{n-1} Q_{n-1} \}$. Only one element from this class will be represented in the quotient space, so it will produce a singular curve with singularity index $k = n/2$. Analogously, the second component $\ell'_2$ corresponds to the class of edges $\{ P_2 Q_2, P_4 Q_4, \ldots, P_n Q_n \}$ and its singularity index equals $k  = n/2$, too. The third component $\ell'_3$ corresponds to the axis of rotation $\rho_k$ and its singularity index equals $k = n/2$.
Through the Heegaard diagram of $M_{25}(n) / \rho_k$ we understand the Heegaard diagram of $M_{25}(2)$ with information about the singular set $\mathcal L_{even}$ presented.
The final step of the proof is based on the fact, proved in Prop.~\ref{prop-M25-2}, that $M_{25}(2)$ is homeomorphic to  $L_{3,1}$.
\qed

Theorem~3.1 from \cite{Cavicchioli2} states that for $n \geqslant 3$ the manifolds $M_{25}(n)$ are hyperbolic, although for $n > 3$ the authors do not present explicitly any proof of hyperbolicity.  Moreover, they give the following volume formula  $\textrm{\rm vol } M_{25} (n) = (n/3) \cdot (4.686034274\ldots)$. However, the given formula turns out to be wrong  even for small $n > 3$: it would be right if the polyhedron $\mathcal P_n$ could be obtained by gluing isometrically  $n$ copies of the $\frac{1}{3}$--piece of the hyperbolic $2 \pi / 3$--icosahedron $\mathcal P_3$.  This is obviously not true, since the dihedral angle  around the image of the axis of rotation $\rho_n$ in $\mathcal P_n / \rho_n$ must be equal to $2 \pi / n$, that is, it differs from $2 \pi / 3$ if $n > 3$. The correct values of $\textrm{\rm vol } M_{25}(n)$ for small  $n$ are presented below.

\begin{prop}
The hyperbolic volumes and the first homology groups of $M_{25}(n)$, for $n\leq 6$, are as follows:
$$
\begin{tabular}{|c|c|c|}
\hline {\rm manifold} & {\rm volume} & {\rm homology group} \cr
\hline  $M_{25}(3)$ & {\rm 4.686034273803\ldots} & $\mathbb Z_2
\oplus \mathbb Z_{18}$ \cr \hline $M_{25}(4)$ & {\rm
3.970289623891\ldots} & $\mathbb Z_3 \oplus \mathbb Z_3 \oplus
\mathbb Z_6$ \cr \hline $M_{25}(5)$ & {\rm 14.319926985892\ldots} &
$\mathbb Z_5 \oplus \mathbb Z_5 \oplus \mathbb Z_{15}$ \cr  \hline
$M_{25}(6)$ & {\rm 14.004768920617\ldots} & $\mathbb Z_8 \oplus
\mathbb Z_{72}$ \cr \hline
\end{tabular}
$$
\end{prop}

\dimo
The results are obtained by using the computer program \emph{Recognizer}.
\qed

\section*{Acknowledgements}

Work performed under the auspices of the GNSAGA of the CNR (National Research Council) of Italy and the Russian Foundation for Basic Research (grants 10-01-00642 and 10-01-91056). The third named author thanks  University of Modena and Reggio Emilia and Abdus Salam School of Mathematical Sciences, GC University Lahore for the hospitality.

\bigskip

\noindent Department of Mathematics, University of Modena and Reggio Emilia, Italy

\noindent Novosibirsk State University, Novosibirsk, 630090, Russia

\noindent Sobolev Institute of Mathematics, Novosibirsk, 630090, Russia

\end{document}